\documentclass[12pt]{amsart}

\usepackage{amssymb}

\usepackage{epsfig,color}

\usepackage{graphicx}

\newtheorem{thm}{Theorem}
\newtheorem{prop}{Proposition}
\newtheorem{defin}{Definition}

\newtheorem{corol}{Corollary}
\newtheorem{cor}{Corollary}
\newtheorem{lem}{Lemma}
\newtheorem{rem}{Remark}

\newtheorem{exa}{Example}

\newcommand{\A}{{\mathcal A}}

\renewcommand{\c}{{\mathcal C}}

\newcommand{\ddd}{{\mathcal D}}

\newcommand{\E}{{\mathcal E}}

\newcommand{\FF}{{\mathcal F}}

\newcommand{\II}{{\mathcal I}}

\newcommand{\N}{{\mathbb N}}\newcommand{\NN}{{\mathbb N}^2}

\newcommand{\QQ}{{\mathbb Q}}

\newcommand{\R}{{\mathbb R}}\newcommand{\RR}{{\mathbb R}^2}

\newcommand{\TT}{{\mathbb T}^2}\newcommand{\ttt}{{\mathcal T}}

\newcommand{\T}{{\mathbb T}}

\newcommand{\Z}{{\mathbb Z}}\newcommand{\ZZ}{{\mathbb Z}^2}

\newcommand{\bo}{\partial} 

\newcommand{\cons}{\mbox{cons}} 
\newcommand{\const}{\mbox{const}} 

\newcommand{\diss}{\mbox{diss}} 

\newcommand{\erg}{\text{erg}}    

\newcommand{\id}{\text{Id}}   
\newcommand{\inter}{\text{interior}}     
\newcommand{\irr}{\text{irr}}  

\newcommand{\leb}{\text{Leb}}     

\newcommand{\rat}{\text{rat}}  

\newcommand{\bg}{\bar{g}}\newcommand{\bG}{\bar{G}}
\newcommand{\bh}{\bar{h}}

\newcommand{\br}{\bar{r}}

\newcommand{\bmu}{\bar{\mu}}

\newcommand{\bPhi}{\bar{\Phi}}

\newcommand{\tx}{\tilde{x}}

\newcommand{\tz}{\tilde{z}}

\newcommand{\tZ}{\tilde{Z}}

\newcommand{\tC}{\tilde{C}}
\newcommand{\tD}{\tilde{D}}
\newcommand{\tE}{\tilde{E}}

\newcommand{\tO}{\tilde{O}}
\newcommand{\tP}{\tilde{P}}

\newcommand{\tR}{\tilde{R}}
\newcommand{\tS}{\tilde{S}}
\newcommand{\tT}{\tilde{T}}
\newcommand{\tX}{\tilde{X}}
\newcommand{\tY}{\tilde{Y}}
\newcommand{\tv}{\tilde{v}}
\newcommand{\vv}{\vec{v}}

\newcommand{\tla}{\tilde{\lambda}}
\newcommand{\tmu}{\tilde{\mu}}
\newcommand{\tnu}{\tilde{\nu}}
\newcommand{\ttau}{\tilde{\tau}}

\newcommand{\tsi}{\tilde{\si}}

\newcommand{\al}{\alpha}
\newcommand{\be}{\beta}
\newcommand{\ga}{\gamma}\newcommand{\Ga}{\Gamma}
\newcommand{\de}{\delta}

\newcommand{\si}{\sigma}\newcommand{\Si}{\Sigma}
\newcommand{\vp}{\varphi}
\newcommand{\het}{\theta}

\begin{document}

\bibliographystyle{plain}

\title[Recurrence and ergodicity]
{On recurrence and ergodicity for geodesic flows on noncompact
periodic polygonal surfaces}

\author{Jean-Pierre Conze and Eugene Gutkin}

\address{IRMAR, CNRS UMR 6625,  Universit\'e de Rennes 1, Campus de Beaulieu, 35042 Rennes Cedex,
France} \email{conze@univ-rennes1.fr}

\address{Copernicus University, Chopina 12/18, Torun 87-100;
IMPAN, Sniadeckich 8, Warszawa 10, Poland}
\email{gutkin@mat.umk.pl,gutkin@impan.pl}



\date{\today}

\begin{abstract}
We study the recurrence and ergodicity for the billiard on
noncompact polygonal surfaces with a free, cocompact action of $\Z$
or $\Z^2$. In the $\Z$-periodic case, we establish criteria for
recurrence. In the more difficult $\Z^2$-periodic case, we establish
some general results. For a particular family of $\Z^2$-periodic
polygonal surfaces, known in the physics literature as the wind-tree
model, assuming certain restrictions of geometric nature,
we obtain the ergodic decomposition of directional billiard
dynamics for a dense, countable set of  directions. This is a
consequence of our results on the ergodicity of $\ZZ$-valued
cocycles over irrational rotations.
\end{abstract}

\maketitle

\tableofcontents

\section*{Introduction}       \label{intro}

Beginning with Boltzmann's {\em ergodic hypothesis}, mathematicians
have been investigating the ergodicity of dynamical systems of
physical origin. Among them are the geodesic flows on riemannian
configuration spaces describing mathematically the physical models
at hand. If the configuration space has a boundary, we arrive at a
billiard.

It is notoriously difficult to study the ergodicity of these
dynamical systems, in particular the famous Boltzmann-Sinai model.
On the contrary, the conservativeness of a dynamical system of this
kind is guaranteed by the Poincar\'e recurrence theorem, provided
its phase space has finite volume. This holds, for instance, if the
configuration space is compact.

The situation changes drastically if the configuration space has
infinite volume. This happens, in particular, if the space is
invariant under a free action of an infinite group, say $\Z^d$. Not
only the ergodicity, but even the conservativeness of these
dynamical systems is a challenging question; it is open in many
relevant examples.

Some physical models correspond to the geodesic flows on polygonal
surfaces invariant under free actions of infinite groups \cite{Gut09}.
We will speak of {\em periodic polygonal surfaces} or {\em $G$-periodic polygonal surfaces},
where $G$ is the group in question.  For instance, the
space of the classical wind-tree model in statistical physics \cite{HW80}
is a $\ZZ$-periodic polygonal surface.

In this work we study the recurrence and ergodicity for  geodesic
flows on noncompact polygons and noncompact polygonal surfaces. For
reader's convenience, we will briefly survey the relevant material
in the compact case. We refer to \cite{Gut84,Gut96} for details. Let
$P$ be a compact polygonal surface, e. g., a polygon. If $P$ is
rational, the study of the geodesic flow on $P$ is equivalent to the
study of the geodesic flow on a compact translation surface, say
$S$. That flow decomposes as a one-parameter family of directional
translation flows, say $T_{\het}^t:S\to S,\,0\le\het\le 2\pi$. The
celebrated result in this subject says that for Lebesgue almost all
directions $\het$ the flows $T_{\het}^t$ are (uniquely) ergodic
\cite{KMS}. This theorem has far reaching applications to the
billiard in irrational polygons \cite{KMS,Vo97}. See
\cite{Gut96,Gut03} for details.

Let now $S$ be a noncompact translation surface. Let
$T_{\het}^t:S\to S,\,0\le\het\le 2\pi,$ be the one-parameter family
of directional translation flows \cite{Gut09}. 
For the purposes of this discussion we assume that $S$ is a
periodic translation surface. The following question naturally
arises: Is there an analog of the unique ergodicity theorem in \cite{KMS} for
noncompact, periodic translation surfaces? This question is mainly
open. In fact, it is not known whether the flows $T_{\het}^t:S\to S$
are conservative for typical directions $\het$. The examples from
our sections~~\ref{one_periodic}, ~~\ref{two_periodic}
and~~\ref{rect_lorenz} might be useful to formulate conjectures
about the flows $T_{\het}^t:S\to S$ for noncompact translation surfaces.

\medskip

We will now informally describe some of our results. Let $\tP$ be a
$\Z$-periodic polygonal surface with a boundary. Let $P=\tP/\Z$ be
the compact quotient. Suppose that the billiard flow on $P$ is
ergodic. Then the billiard flow on $\tP$ is conservative. See
Theorem~~\ref{one_period_thm}.

With any polygon $O$ inside the unit square we associate the
$\Z$-periodic strip $\tP_{O}$ with polygonal obstacles. Then for a
dense $G_{\de}$-set of obstacles, the billiard flow on $\tP_{O}$ is
conservative. See Theorem~~\ref{dens_G_del_thm}. Let $O$ be an {\em
irrational polygon}. Suppose that its angles admit a
superexponentially fast approximation by numbers in $\pi\QQ$. Then
the geodesic flow on $\tP_O$ is conservative. See
Theorem~~\ref{vorob_cond_thm}.

Let $\tP=\tP(a,b)$ be the {\em rectangular Lorenz gas} obtained by
deleting from $\RR$ the $\ZZ$-periodic family of $a\times b$
rectangles. This corresponds to the {\em wind-tree model}
\cite{HW80}. Let $p,q\in\N$ be relatively prime; denote by
$\tT(a,b;p,q)$ the billiard flow in the direction $\arctan{(q/p)}$. We
say that the {\em obstacles are small} if $qa+pb \le 1$. Assuming
the small obstacles condition, we analyze the flow $\tT(a,b;p,q)$.
Let $\tT_{\cons}(a,b;p,q)$ and $\tT_{\diss}(a,b;p,q)$ be the
conservative and the dissipative parts of $\tT(a,b;p,q)$
respectively. We show that $\tT_{\diss}(a,b;p,q)$ is trivial  iff
$qa+pb=1$. Assume now that $a/b$ is irrational. Then 
 we obtain an ergodic decomposition of $\tT_{\cons}(a,b;p,q)$. 
The $2pq$ ergodic components are isomorphic; they have a simple geometric meaning. 
Thus, $\tT_{\cons}(a,b;p,q)$ is a
finite multiple of an ergodic flow. See
Theorem~~\ref{erg_decom_thm}, Theorem~~\ref{bil_erg_decom_thm}, and
Proposition~~\ref{bil_erg_decom_cor}.

\begin{figure}[htbp]
\begin{center}
\includegraphics[scale=0.5]{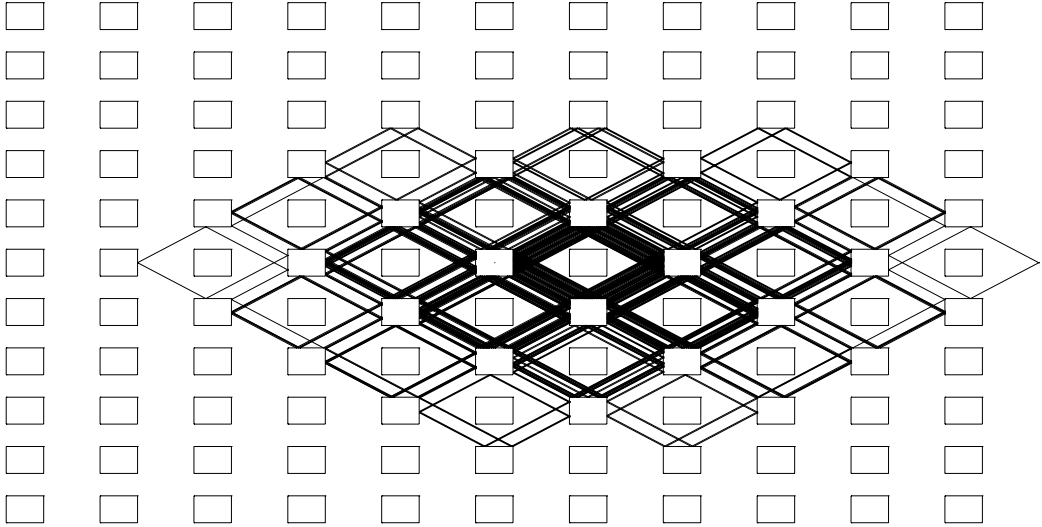}
\end{center}
\vspace{4.5 cm}
\caption{\it An orbit for the flow $\tT_{\cons}(a,b;1,1)$.}
\label{fig7}
\end{figure}

For instance, the conservative part of the wind-tree billiard flow in direction $\pi/4$
is the flow $\tT_{\cons}(a,b;1,1)$; it
has two ergodic components. Figure~~\ref{fig7} shows a typical orbit
of this flow. It encounters only a half of the set of rectangular 
obstacles. Loosely speaking, the orbit skips every other obstacle. The skipped
obstacles are visited by a typical orbit from the other ergodic component of
the flow. 

\newpage

This is a special case of the general situation, as we explain in 
Theorem~~\ref{erg_decom_thm} and Theorem~~\ref{bil_erg_decom_thm}.

\vspace{1.5 cm}

We will now comment on our methods. The billiard flows are
suspensions of billiard maps. Exploiting the periodicity under a group $G$, 
we identify a billiard map of this kind with a skew product over an interval
exchange with the fiber $G$. In our
setting, we obtain skew products with the fibers $\Z$ and $\ZZ$.
Thus, we reduce the questions concerning the ergodicity of the billiard
on noncompact polygonal surfaces to the ergodicity of particular
$\Z^d$-valued cocycles over interval exchanges.

This work is not concerned with the subject of ergodicity for cocycles over
general interval exchanges. We establish the ergodicity of a class
of $\Z^2$-valued cocycles over irrational rotations. See
Theorems~~\ref{wdd-erg-thm},~~\ref{erg-beta}, and ~~\ref{ergo-Psi}
in  section~~\ref{ergodic}. Theorem~~\ref{wdd-erg-thm} states that
a cocycle of this kind is ergodic if the continued fraction
decomposition of the rotation number satisfies certain genericity
assumptions which hold for almost all numbers.
Theorem~~\ref{erg-beta} and Theorem~~\ref{ergo-Psi} strengthen
Theorem~~\ref{wdd-erg-thm} by removing these assumptions.

We apply these results in section~~\ref{rect_lorenz} to obtain the
ergodic decompositions of directional flows in the wind-tree model. 
Let $T^t$ be one of the $2pq$ geometric components of $\tT_{\cons}(a,b;p,q)$.
The Poincar\'e map for $T^t$ is a skew product with fibre $\ZZ$ over a circle
rotation; this rotation is irrational iff $a/b$ is irrational. 
The corresponding $\ZZ$-valued cocycle
belongs to the class of cocycles studied in section~~\ref{ergodic}.
Theorem~~\ref{ergo-Psi} implies the ergodicity of $T^t$.

\medskip

We will now outline the structure of our exposition.
Section~~\ref{cadre} contains the information to be used in the body
of the paper. In section~~\ref{ergodic_sub} we review the material
on skew products, cocycles, recurrence and transience. In
section~~\ref{trans_surf_sub} we establish the framework of
polygonal surfaces. In section~~\ref{bill_sub} we recall the basic
facts about the billiard flow and the billiard map.

Section~~\ref{one_periodic} is about the $\Z$-periodic case. This
framework is naturally divided into two extreme situations: the
generic case and the rational case. In the former situation we show
the conservativeness and some ergodic properties for the generic
$\Z$-periodic polygonal surface; in the latter we establish these
properties for the directional flows in almost all directions.

Section~~\ref{two_periodic} and section~~\ref{rect_lorenz} are
devoted to the billiard on $\Z^2$-periodic polygonal surfaces. In
section~~\ref{two_periodic} we establish some ergodic properties for
arbitrary $\Z^2$-periodic polygonal surfaces. In
section~~\ref{rect_lorenz} we consider a particular family of such
surfaces: The rectangular Lorenz gas or the wind-tree model of the
physics literature. Applying the results of section~~\ref{ergodic},
we obtain the ergodic decomposition of directional billiard flows
for a dense, countable set of directions, under restrictions
of geometric nature, namely the smallness of obstacles condition.
Section~~\ref{ergodic} is a study of ergodicity for a class of
cocycles over irrational rotations. The results, besides being of
interest on their own, are instrumental for the material in
section~~\ref{rect_lorenz}.

\section{The setting and preliminaries}    \label{cadre}
For convenience of the reader, we recall the basic material about
recurrence \cite{Aa97, CorFomSin, Sc77}. We will consider two kinds
of {\em dynamical systems}: transformations and flows. In the former
case, we have the standard Borel space $(X, {\mathcal A})$ endowed
with a possibly infinite measure $\nu$, and a transformation $\tau:
(X, {\mathcal A})\to(X, {\mathcal A})$ preserving $\nu$. For
simplicity, we will assume that $\tau$ is invertible. The setting
for flows is analogous \cite{CorFomSin}. We will use the notation
$(X,\tau,\nu)$ (resp. $(Y,T^t,\mu)$) for transformations (resp.
flows).

The dynamical system $(X,\tau,\nu)$ is {\em recurrent} or {\em
conservative} if for every measurable set $B\subset X$ and for
$\nu$-a.e. point $x\in B$ there is $n=n(x)>0$ such that $\tau^n x
\in B$. Recurrence for flows $(Y,T^t,\mu)$ is defined analogously. A
dynamical system uniquely decomposes as a disjoint union of the {\em
conservative} part and the {\em dissipative part}. For simplicity,
we describe this decomposition only for a transformation,
$(X,\tau,\nu)$. The conservative, dissipative subsets $C,D\subset X$
are measurable and $\tau$-invariant. If $\nu(C)>0$ then
$(C,\tau,\nu)\subset(X,\tau,\nu)$ is recurrent. Suppose that
$\nu(D)>0$. Then there is a measurable set $A\subset D$ such that $D
= \cup_{n \in \Z} \tau^n A$; moreover, $\tau^p A \cap \tau^q A =
\emptyset$ for $p \ne q$.

If $\nu(X)=\infty$, this decomposition is, in general, nontrivial.
We will use the following observation. Let $(Y,T^t,\mu)$ be a flow,
let $X\subset Y$ be a {\em cross-section}, and let $(X,\tau,\nu)$ be
the induced transformation. Then the transformation $(X,\tau,\nu)$
is recurrent iff the flow $(Y,T^t,\mu)$ is recurrent.

The geometric spaces that we work with in the body of the paper are
differentiable manifolds, possibly with boundary and corners. The
transformations and flows are piecewise differentiable.

\subsection{Ergodic theory for skew products}
\label{ergodic_sub}
\hfill \break We will represent our dynamical systems as skew
products over dynamical systems with finite invariant measures.
Their fibers will be infinite abelian groups.

In this section we recall the relevant material on $G$-valued
cocycles, where $G$ is an infinite abelian group; we will write the
group operation additively. We restrict the discussion mostly to the
groups $G=\R^m\times\Z^n$. We denote by $\leb_G$ a Haar measure on
$G$, suppressing the subscript if the group is clear from the
context.

\begin{defin}   \label{cocycle_def}
Let $(X,\tau,\nu)$ be a dynamical system, and let $\vp:X\to G$ be a
measurable function. It determines a {\em cocycle} $\vp(n, x)$, also
denoted by $\vp_n(x)$ or simply $(\vp_n)$, as follows. We set
$\vp(0, x) = 0$. For $n\ne 0$ we set
\begin{equation}    \label{cocycle_eq}
\vp_n(x) = \sum_{j=0}^{n-1} \vp(\tau^j x), \mbox{if}\ n > 0; \
\vp_n(x) = - \sum_{j=n}^{-1} \vp(\tau^j x), \mbox{if}\ n < 0.
\end{equation}
\end{defin}


Thus, $\vp(n, x)$ are the ergodic (or Birkhoff) sums of $\vp$ with
respect to the transformation $(X,\tau,\nu)$. The cocycle $(\vp_n)$
can be viewed as the random walk on the group $G$ driven by the
dynamical system $(X,\tau,\nu)$. Set
\begin{equation}     \label{produit_gauche_eq}
\ttau(x,g) = (\tau x, g+\vp(x)).
\end{equation}
Then $\ttau$, or $\tau_\vp$ to emphasize the dependence on $\vp$, is
a transformation of $\tX=X\times G$ preserving the product measure
$\tnu =\nu\times\leb$. The dynamical system $(\tX,\ttau,\tnu)$ (or
$(\tX,\tau_\vp,\tnu)$) is the {\it skew product} over $(X,\tau,\nu)$
with the {\it displacement function} $\vp$.
Equation~~\eqref{cocycle_eq} corresponds to the iterates of $\ttau$.
Namely, for $n\in\Z$ we have
\begin{equation}     \label{prod_gauche_iter_eq}
\ttau^n(x,g)=(\tau^n(x),g+\vp_n(x)).
\end{equation}

\begin{defin}    \label{recurr1_def}
Let $(X,\tau,\nu)$ be a dynamical system with $\nu(X)<\infty$. Let
$\vp:X\to G$ be a measurable function. The cocycle $(\vp_n)$ is {\it
transient at $x\in X$} if $\vp_n(x)\to \infty$; otherwise, the
cocycle is {\em recurrent at $x$}.\footnote{We will also say that
$x$ is a transient (resp. recurrent) point for the cocycle.} The
cocycle is recurrent (resp. transient) if it is recurrent (resp.
transient) at a.e. $x \in X$.
\end{defin}

We point out a subtlety in Definition~~\ref{recurr1_def}. Let
$\al_n, \be_n$ be $\Z$-valued cocycles over $(X,\tau,\nu)$. Their
direct sum $\vp_n=(\al_n,\be_n)$ is a $\ZZ$-valued cocycle. The
recurrence of $\al_n,\be_n$ does not necessarily imply that $\vp_n$
is recurrent. See, e. g., \cite{ChCo09} for an example.

The sets of transient and recurrent points for a cocycle are
measurable and invariant. Hence, any cocycle $\vp_n$ over an ergodic
$(X,\tau,\nu)$ is either recurrent or transient. Let $(X,\tau,\nu)$
be arbitrary, let $(\vp_n)$ be a cocycle, and let $R\subset X$ be
the set of recurrent points for $\vp_n$. Suppose that $\nu(R)>0$,
and let $\nu_{R}$ be the restriction of $\nu$ to $R$; set
$\tnu_{R}=\nu_{R}\times\leb$, $\tR=R\times G$. Then the skew product
$(\tR,\ttau,\tnu_{R})$ is a conservative dynamical system
\cite{Sc77}. Assume, moreover, that $X$ is a separable metric space
and that $\tau:X\to X$ is compatible with the topological
structure.\footnote{These assumptions will be satisfied in our
applications.} Then for a.e. $x\in R$ there is an infinite sequence
$n_k=n_k(x)$ such that $\tau^{n_k} x \to x$ and $\vp(n_k, x) \to 0$
\cite{Sc77}. Therefore, for almost every point in $R$ the sequence
$(\vp_k(x))_{k \ge 0}$ visits arbitrarily close to $0\in G$. If
$G=\Z^d$, then we have $\vp_k(x) = 0$ infinitely many times. Suppose
now that $(\tX,\ttau,\tnu)$ comes from a billiard model. Then $R$ is
the conservative part of the billiard phase space. The billiard ball
emanating from a point in $R$ almost surely returns infinitely often
to the obstacle from which it started, and arbitrarily close to the
point of departure.

\vskip 3mm We refer the reader to \cite{Sc77} for criteria of
recurrence for cocycles.  See also \cite{Ke75}, \cite{At76}. The
following lemma from \cite{ChCo09} gives a simple
sufficient condition for recurrence of $\R^d$-valued cocycles.

\begin{lem}       \label{recurrence_lem}
Let $(X,\tau,\nu)$ be a dynamical system with a finite measure, let
$\vp:X\to\R^d$ be a measurable function, and let $\vp(n,x)$ be the
corresponding cocycle.

Let $| \ \ |$ be a norm on $\R^d$. Suppose that there exists a
strictly increasing sequence of integers $k_n$ and a sequence of
nonnegative functions $\delta_n(x)$ converging to $0$ for almost
every $x$ such that
\begin{equation}    \label{lyapunov_eq}
\lim_n \nu(\{x: |\varphi(k_n, x)| \ge \delta_n(x)
n^{\frac{1}{d}}\})=0.
\end{equation}

Then the cocycle is recurrent. \end{lem}

Let $(X,\tau,\nu)$ be a dynamical system, and let $\vp:X\to G$ be a
measurable function. Recall that $\vp$ is a {\em coboundary} if
there exists a measurable function $\psi:X\to G$ such that
$\vp=\psi-\tau\psi$. Let $(\tX,\ttau,\tnu)$ be the dynamical system
defined by equation~~\eqref{produit_gauche_eq}. Let $H\subset G$ be
a closed subgroup; set $\tX_H=X\times
G/H,\tnu_H=\nu\times\leb_{G/H}$. Define $\ttau_ H: \tX_H \to \tX_H$
by $\ttau_ H(x,g +H) = (\tau x, (\vp(x)+ g) + H)$. Then
$(\tX_H,\ttau_ H,\tnu_H)$ is the skew product over $(X,\tau,\nu)$
with the fibre $G/H$ and the displacement function
$\vp_H(x)=\vp(x)+H$.

\begin{lem}             \label{cobord_lem}
If there is a closed, proper subgroup $H\subset G$ such that the
dynamical system $(\tX_H,\ttau_H,\tnu_H)$ is ergodic, then $\vp$ is
not a coboundary.
\begin{proof}
Assume the opposite, and let $\vp=\psi-\tau\psi$ where $\psi:X\to G$
is a measurable function. Set $\Psi(x,g)=(x,\psi(x)+g)$. Then
$\Psi:X\times G\to X\times G$ is an automorphism of the measure
space $(X\times G,\tnu)$; it conjugates $\ttau$ and the product
transformation $\tau\times\id$. Dividing by $H$, we obtain the
automorphism $\Psi_H:X\times G/H\to X\times G/H$ conjugating
$\ttau_H$ and $\tau\times\id_{G/H}$. This contradicts the ergodicity
of $\ttau_H$.
\end{proof}
\end{lem}

We introduce a terminology to express the property that a walk in
$G$ visits any compact set with zero asymptotic frequency.

\begin{defin}   \label{NRUO_def}
Let $(X,\tau,\nu)$ be a dynamical system with finite invariant
measure, let $G\subset\R^d$ or $\Z^d$ be a closed subgroup, and let
$\vp:X\to G$ be a measurable function. The associated cocycle
$(\vp_n)$ is {\em zero-recurrent} if it is recurrent, and for a.e.
point $x \in X$ and any compact set $K\subset G$ we have
\begin{equation}     \label{freq0}
\lim _{n\to\infty}\left\{{1\over n} \sum_{k = 0}^{n-1} 1_K(\vp(k,
x))\right\} = 0.
\end{equation}
\end{defin}

\begin{prop}  \label{recurrence_prop}
Let $(X,\tau,\nu)$ be a dynamical system with finite invariant
measure, let $\vp:X\to\R$ be an integrable function, let $(\vp_n)$
be the associated cocycle, and let $(\tX,\ttau,\tnu)$ be the
corresponding skew product. Denote by ${\mathcal J}$ the
$\sigma$-algebra of measurable $\tau$-invariant subsets. Let
${\mathbb E}(\vp |{\mathcal J})$ be the conditional expectation of
$\vp$ with respect to ${\mathcal J}$. Let $R\subset X$ be the set of
recurrent points for the cocycle $(\vp_n)$.

Then the following properties hold.

\begin{itemize}

\item[1.] The set $R$
and the set $\{x: {\mathbb E}(\vp |{\mathcal J})(x) = 0 \}$ coincide
up to a set of $\nu$-measure zero.

\item[2.] If the dynamical system $(X,\tau,\nu)$ is ergodic and
$\int_X \vp \ d\nu = 0$, then the cocycle $(\vp_n)$ is recurrent.

\item[3.]  If, moreover, $\vp$ is not a coboundary, then the
cocycle $(\vp_n)$ is zero-recurrent.
\end{itemize}
\begin{proof}
Let $A,B\subset X$ be measurable sets. By $A=B$ we will mean that
$A$ and $B$ are equal in the measure-theoretic sense.\footnote{I.
e., their symmetric difference is of $\nu$-measure $0$.} We will
also say, simply, that $A$ and $B$ coincide.

The first two claims are classical. For $\nu$-a.e. $x$ the
conditional expectation ${\mathbb E}(\vp |{\mathcal J})(x)$ is
defined, and, by the Birkhoff ergodic theorem,
$\frac{1}{n}\varphi(n,x) \to {\mathbb E}(\vp |{\mathcal J})(x)$.
Hence, the $\tau$-invariant sets $\{x:{1\over n} \varphi(n,x)
\rightarrow 0\}$ and  $\{x: {\mathbb E}(\vp |{\mathcal J})(x) = 0\}$
coincide. We denote this $\tau$-invariant set by $X_0$. By
Lemma~~\ref{recurrence_lem}, the cocycle $(\vp_n)$ is recurrent on
$X_0$. By the Birkhoff ergodic theorem, almost every $x\in
X\setminus X_0$ is transient for $(\vp_n)$; moreover, on $X\setminus
X_0$, the cocycle has linear dissipation. Thus, $R=X_0=\{x: {\mathbb
E}(\vp |{\mathcal J})(x) = 0\}$, proving claim 1. Claim 2 directly
follows from claim 1. We will now prove claim 3.

Set $\tX=X\times\R$ and let $(\tX,\ttau,\tnu)$ be the skew product
equation~~\eqref{produit_gauche_eq}. Let $K\subset \R$ be any
compact. Set
$$u_K(x,g)=\lim_{n\to\infty} {1\over n} \sum_{k = 0}^{n-1} 1_K(\vp(k,
x) + g).$$

The ergodic theorem applied to $(\tX,\ttau,\tnu)$ ensures the
existence of the limit for a.e. $(x,g)$ and that $u_K$ is an
integrable, nonnegative, $\ttau$-invariant function on $\tX$.
Suppose that $u_K\ne 0$ on a set of positive measure. Then $u_K\tnu$
is a finite $\ttau$-invariant measure on $\tX$, absolutely
continuous with respect to $\tnu$. Since $(X,\tau,\nu)$ is ergodic,
this implies that $\varphi$ is a coboundary \cite{Co79}, contrary to
the assumption. Thus, for any compact $K\subset \R$ we have $\lim
\frac1n \sum_{k = 0}^{n-1}1_K(\vp(k, x) + g) = 0$ for $\tnu$-a.e.
$(x,g)\in\tX$. Since $\R$ is a countable union of compacta, there
exists $\tY\subset\tX$, $\tnu(\tX\setminus\tY)=0$, such that for
$(x,g)\in\tY$ and any compact $K\subset\R$ we have $u_K(x,g)=0$.

By Fubini's theorem, there exists $g_0\in\R$ such that, for $\nu$
a.e. $x\in X$ we have $\lim \frac1n \sum_{k = 0}^{n-1} 1_K(\vp(k, x)
+ g_0) = 0$ for any compact $K\subset\R$. But $g\mapsto g+g_0$ is a
self-homeomorphism of $\R$.
\end{proof}
\end{prop}

\begin{rem}   \label{zero_rec_rem}
{\em The same argument proves claim 3 for $\R^d$-valued cocycles.
See \cite{Co79} for a generalization to cocycles with values in
locally compact groups.

}
\end{rem}

Proposition~~\ref{recurrence_prop} and Lemma~~\ref{cobord_lem} imply
the following.

\begin{corol}
Let $(X,\tau,\nu)$ be a dynamical system with finite measure. Let
$\vp:X\to\R$ be a measurable function such that $\int_X \vp\, d\nu =
0$, and let $(\varphi_n)$ be the associated cocycle. Let
$(\tX,\ttau,\tnu)$ be the skew product corresponding to $\vp$. For
$p\in\R$ let $(\tX_p,\ttau_p,\tnu_p)$ be the reduction of
$(\tX,\ttau,\tnu)$ with respect to the subgroup $H=p\Z$, as in
Lemma~~\ref{cobord_lem}.

Suppose that for some $p\ne0$ the dynamical system
$(\tX_p,\ttau_p,\tnu_p)$ is ergodic. Then the cocycle $(\varphi_n)$
is zero-recurrent.
\end{corol}

We will use the following proposition.
\begin{prop} \label{lem-k_n}
Let $(X,\tau,\nu)$ be an ergodic dynamical system with a finite
measure, let $\vp:X\to\R^d$ be a measurable function, and let
$\vp(n,x)$ be the corresponding cocycle. Then the following
dichotomy holds: i) The cocycle $(\varphi_n)$ is recurrent; ii) Let
$k_n\in\N$ be any strictly increasing sequence. Then there exists
$c>0$\footnote{In general, it depends on the sequence.} such that
for a.e. $x\in X$ we have
\begin{eqnarray}   \label{lim-k_n}
\limsup_n (n^{-1/d} |\varphi(k_n, x)|) = c.
\end{eqnarray}
\begin{proof}
The quantity $\limsup_n (n^{-{1\over d}} |\varphi(k_n, x)|)$ is an
invariant, measurable function. By ergodicity, it is equal to a
constant $c \ge 0$. The claim is now immediate from
Lemma~~\ref{recurrence_lem}.
\end{proof}
\end{prop}

\subsection{Noncompact, periodic polygonal surfaces} \label{trans_surf_sub}
\hfill \break We will now establish the geometric framework for our
study. Let $G$ be an infinitely countable group acting freely and
cocompactly by isometries on a noncompact riemannian manifold
$\tP$.\footnote{In general, with boundary and corners.} Then
$P=\tP/G$ is compact; the projection $p:\tP\to P$ is a riemannian
covering. Let $U\tP,UP$ be the unit tangent bundles for $\tP,P$; let
$\tT^t,T^t$ be the respective geodesic flows; let $\tmu,\mu$ be the
Liouville measures for $U\tP,UP$ respectively. The action of $G$ on
$\tP$ uniquely extends to a free, cocompact action on $U\tP$. We
have $UP=U\tP/G$; let $q:U\tP\to UP$ be the projection. Then
$q:(U\tP,\tT^t,\tmu)\to(UP,T^t,\mu)$ is a covering of flows.

Let $X\subset UP$ be a compact submanifold which is a cross-section
for $(UP,T^t,\mu)$. Then  the manifold $\tX=p^{-1}(X)\subset U\tP$
is a cross-section for the flow $(U\tP,\tT^t,\tmu)$. Let $\nu,\tnu$
be the induced measures on $X,\tX$ respectively; let $\tau:X\to
X,\,\ttau:\tX\to \tX$ be the respective Poincar\'e maps. Then $\tX =
X\times G$ measure theoretically, and $\tnu=\nu\times\leb$. There is
a unique mapping $\vp:X\to G$ such that $(\tX,\ttau,\tnu)$ is the
skew product over $(X,\tau,\nu)$ with the displacement function
$\vp$. See equation~~\eqref{produit_gauche_eq}.

Let the manifold $\tP$ be a {\em noncompact polygonal surface}. See
\cite{Gut84,Gut09} for the background. We will say that $\tP$ is a
{\em $G$-periodic} polygonal surface. When the group $G$ is
implicit, we will say that $\tP$ is a periodic polygonal surface.
When $G=\Z$ or $G=\ZZ$, we will say that $\tP$ is $\Z$-periodic or
$\ZZ$-periodic respectively. The projection $p:\tP\to P$ is a {\em
covering of polygonal surfaces} \cite{Gut09}.

If $\tP$ (or, equivalently, $P$) is a {\em rational polygonal
surface}, we associate with it a finite subgroup $\Ga=\Ga(P)\subset
O(2)$. For $N\ge 1$ let $R_N\subset O(2)$ be the dihedral group of
order $2N$, i. e., the group generated by two orthogonal
reflections, with the angle $\pi/N$ between their axes. If $\bo
P\neq\emptyset$, then $\Ga(P)=R_N$, where $N$ is determined by $P$.
Note that $R_1$ consists of a reflection and the identity. The
associated {\em translation surface} $S=S(P)$ \cite{GJ} is a compact
riemann surface endowed with a riemannian metric, flat everywhere
except for a finite number of {\em cone points}. The group $\Ga(P)$
acts on $S(P)$ by isometries, and we have $P=S(P)/\Ga(P)$
\cite{Gut84}. Thus, a polygonal surface is a translation surface iff
$\Ga(P)=\id$. Let $\Ga=\Ga(P)$ and let $S=S(P)$. Then there is a
unique noncompact, $G$-periodic translation surface $\tS$ such that
the following conditions hold. The groups $\Ga$ and $G$ act on $\tS$
by isometries; the two actions commute. We have $\tS/G=S$,
$\tS/\Ga=\tP$, $S/\Ga=P$, $\tP/G=P$; the projections $\tS\to S$,
$\tP\to P$, $\tS\to\tP$, $S\to P$ are compatible.


The group $\Ga$ acts on the unit circle $U\subset\RR$; let $U/\Ga$
be the quotient. The flows $(U\tP,\tT^t,\tmu)$ and $(UP,T^t,\mu)$
decompose as one-parameter families of {\em directional geodesic
flows} $(U\tP_{\het},\tT^t_{\het},\tmu_{\het})$ and
$(UP_{\het},T^t_{\het},\mu_{\het})$ where $\het\in U/\Ga$. The
projection $q$ is compatible with the decompositions, inducing the
directional projections $q_{\het}:(U\tP_{\het},\tT^t_{\het},
\tmu_{\het})\to(UP_{\het},T^t_{\het},\mu_{\het})$. The flow
$(UP_{\het},T^t_{\het},\mu_{\het})$ is naturally isomorphic to the
linear flow in direction $\het$ on the translation surface $S$; the
measure $\mu_{\het}$ corresponds to the Lebesgue measure on $S$.
Analogously, $(U\tP_{\het},\tT^t_{\het},\tmu_{\het})$ is the {\em
linear flow in direction $\het$ on the periodic translation surface}
$\tS$ \cite{Gut09}.

\begin{exa}    \label{strip_exa}
{\em Let $B\subset\RR$ be the horizontal strip bounded by the lines
$\{y=0\}$ and $\{y=1\}$. For $0\le a,b <1,a+b>0$ let $R_0=R_0(a,b)$
be the $(a\times b)$-rectangle centered at $(1/2,1/2)$. Set
$\tP(a,b)=B\setminus\cup_{k\in\Z}(R_0+(k,0))$. The $\Z$-periodic
polygonal surface $\tP(a,b)$ is a strip with a periodic sequence of
rectangular obstacles. See figure~~\ref{rect_loren_band}. Let $Q$ be
the unit square $0\le x,y \le 1$; let $C$ be the cylinder obtained
by identifying the vertical sides of $Q$. Then $P(a,b)=C\setminus
R_0(a,b)$ is the unit cylinder with a rectangular obstacle. In the
limit cases $a=0$ or $b=0$ the obstacles degenerate into {\em
barriers}. We have $\Ga=R_2$ if $b\ne 0$ and $\Ga=R_1$ if $b=0$.

\begin{figure}[htbp]
\begin{center}
\input{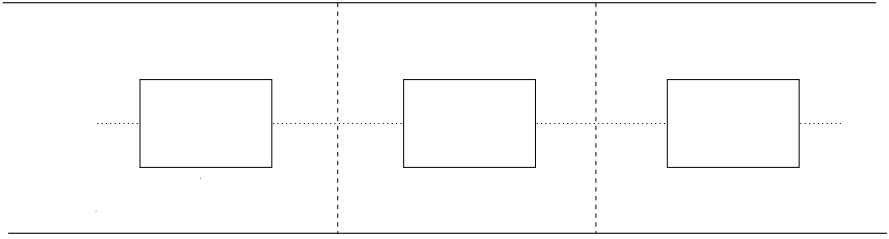}
\caption{\it An infinite band with a periodic configuration of
rectangular obstacles.} \label{rect_loren_band}
\end{center}
\end{figure}

Set $P=P(a,b)$ and $S=S(a,b)$. Let $a,b\ne 0$. The translation
surface $S$ is constructed from 4 copies of $Q\setminus R_0(a,b)$
via identifications of their sides shown in
figure~~\ref{quotient_S(a,b)}. It has $4$ cone points with cone
angles $6\pi$. This yields $g(S)=5$ \cite{GJ}. The genus of $S$ can
also be computed directly from the angles in $P$ \cite{Gut84}. The
same analysis applies when $a=0$ or $b=0$. If $a=0$, the surface $S$
is made from $4$ rectangles with vertical barriers. It has $4$ cone
points with cone angles $4\pi$, yielding $g(S)=3$. If $b=0$, then
$S$ is  made from $2$ rectangles with horizontal barriers. There are
$2$ cone points with cone angles $4\pi$, yielding $g(S)=2$.

Set $\tP=\tP(a,b)$ and $\tS=\tS(a,b)$. The noncompact translation
surface $\tS$ is obtained by analogous identifications of pairs of
sides in the disjoint union of  4 copies of $\tP$. Since $\tP$ is
$\Z$-periodic, and since these identifications are compatible with
the action of $\Z$, the translation surface $\tS$ is $\Z$-periodic.
It has infinite genus. }
\end{exa}

\begin{figure}[htbp]
\begin{center}
\input{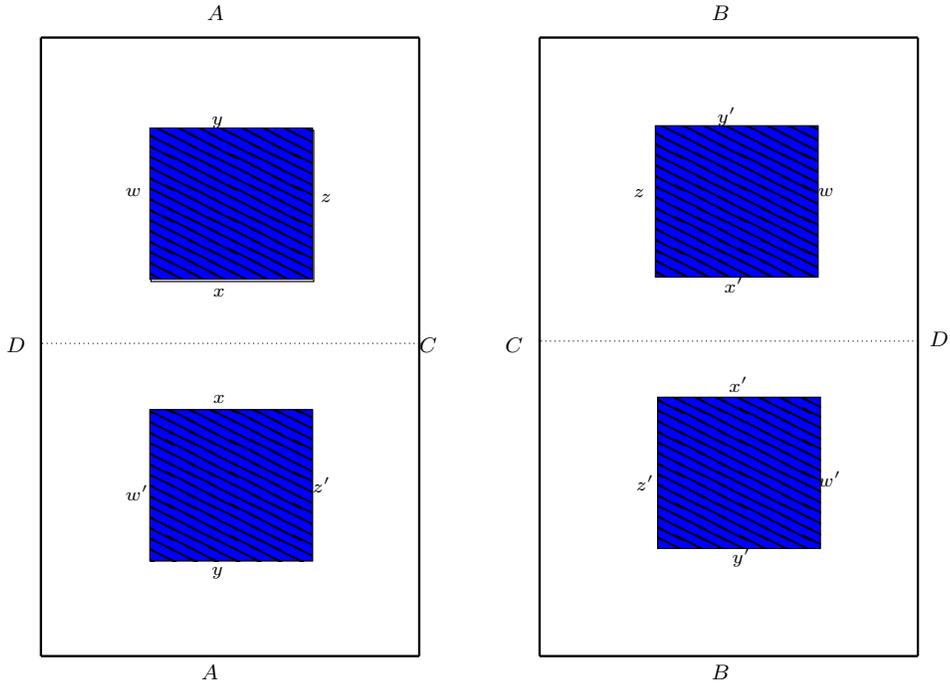}
\caption{\it The translation surface made from {\em 4} copies of
$Q\setminus R_0(a,b)$ by identifying the sides bearing the same
labels.} \label{quotient_S(a,b)}
\end{center}
\end{figure}

\subsection{The billiard flow and
the billiard map} \label{bill_sub}
\hfill \break Let $P$ be a compact polygonal surface, and let $\bo
P$ be its boundary. Orbits of the geodesic flow $(UP,T^t,\mu)$,
viewed as curves in $P$, are the geodesics. We will use the term
{\em billiard curves} for those geodesics that intersect $\bo P$ at
regular points.

Virtually every polygonal surface has singular points. Points
$z\in\inter(P)$ (resp. $z\in\bo(P)$) are singular if they are cone
points (resp. corner points). A geodesic in $P$ may start or end at
a singular point, but it cannot pass through a singular point. A
phase point $v\in UP$ is singular if the geodesic it defines arrives
at a singular point, and hence is defined on a proper subinterval of
$(-\infty,\infty)$. Therefore the geodesic flow
$T^t,-\infty<t<\infty,$ is defined only on the set of regular (i.
e., non-singular) points in $UP$. We will use the terms {\em regular
and singular sets} for the sets of regular and singular points,
respectively.

The singular set has codimension one; thus, the Liouville measure of
the singular set is zero, and the regular set has full measure. As
is usual in the billiard literature \cite{Gut84,GK}, we will use the
notation $(UP,T^t,\mu)$ for the geodesic flow on the regular set.
Analogous notational conventions are used for the billiard map,
which we will now define. The reader should mentally substitute
``regular phase points'' whenever we speak of phase points in what
follows. The billiard map will be defined on the regular set which
has full measure in the canonical cross-section that we will now
describe.

Let $A\subset P$. We denote by $U_AP\subset UP$ the set of vectors
$v\in UP$ whose base points belong to $A$. Let $BUP\subset UP$ be
the smallest $T^t$-invariant set containing $U_{\bo P}P$. If $\bo
P\ne\emptyset$, we set $\mu_B=\mu|_{BUP}$. By definition, $U_{\bo
P}P$ is a cross-section for the flow $(BUP,T^t,\mu_B)$. Let $\nu$ be
the induced measure on $U_{\bo P}P$. We will use the following
terminology: The flow $(BUP,T^t,\mu_B)$ is the {\em billiard flow}
of $P$; the set $U_{\bo P}P\subset BUP$ is the {\em standard
cross-section} for the billiard flow; the induced transformation
$(U_{\bo P}P,\tau,\nu)$ is the {\em billiard map} for $P$.

Let now $P$ be a rational polygonal surface with a boundary, and let
$\Ga(P)=R_N$. We identify $U/\Ga$ and $[0,\pi/N]$. For
$\het\in[0,\pi/N]$ set $U_{\bo P}P_{\het}=UP_{\het}\cap U_{\bo P}P$,
$BUP_{\het}=BUP\cap UP_{\het}$. Suppose that
$\mu_{\het}(BUP_{\het})>0$. Let $b\mu_{\het}$ be the restriction of
$\mu_{\het}$ to $BUP_{\het}$; then
$(BUP_{\het},T^t_{\het},b\mu_{\het})$ is the {\em billiard flow in
direction $\het$}. The set $U_{\bo P}P_{\het}\subset BUP_{\het}$ is
the {\em standard directional cross-section}. Let $\nu_{\het}$ be
the induced measure on $U_{\bo P}P_{\het}$. The induced dynamical
system $(U_{\bo P}P_{\het},\tau_{\het},\nu_{\het})$ is the {\em
directional billiard map} \cite{Gut09}. When $P\subset\RR$ is a
compact polygon, this is the standard terminology \cite{Gut03}.

Let $S$ be a compact translation surface. Let $O\subset S$ be a
polygon, not necessarily connected, such that $S\setminus O$ is
connected. The polygonal surface $P=S\setminus\inter(O)$ is a {\em
translation surface with polygonal obstacles}. Some of the
components of $O$ may be linear segments, hence we will also speak
of {\em translation surface with polygonal obstacles and/or
barriers}.

\begin{lem}          \label{bil_type_lem}
1. Let $P$ be a compact translation surface with polygonal obstacles
and/or barriers. Then $\mu(UP\setminus BUP)=0$. 2. Let $P$ be a
compact, rational polygonal surface with a boundary; let $R_N$ be
the corresponding reflection group. Let
$\E_{\erg}(P)\subset[0,\pi/N]$ be the set of uniquely ergodic
directions. i) If $\bo P$ contains intervals with distinct
directions, then for every $\het\in\E_{\erg}(P)$ we have
$BUP_{\het}=UP_{\het}$. ii) Suppose that $\bo P$ consists of
intervals with the same direction, say $\het_0$. Then for
$\het\in\E_{\erg}(P)\setminus\{\het_0\}$ we have
$BUP_{\het}=UP_{\het}$.
\begin{proof}
1. Let $P=S\setminus O$. Let $\ga$ be an infinite geodesic in $P$
that does not intersect $O$. Then $\ga$ is an infinite geodesic in
$S$, and $\ga\cap O=\emptyset$. If $\inter(O)\ne\emptyset$, then
$\ga$ is not dense, and hence its direction is not minimal. Suppose
$\inter(O)=\emptyset$, i. e., $O$ consists of barriers. Assume that
$O$ contains segments of distinct directions. A geometric argument
which we leave to the reader implies that $\ga$ is not dense in $S$,
and hence its direction is not minimal. Let $O$ consist of segments
with direction $\het_0$. Then the direction of $\ga$ is either
nonminimal or it is $\het_0$. But the set of nonminimal directions
is countable, implying the claim.

2. There is a compact translation surface with obstacles and/or
barriers, say $S\setminus O$, and a finite covering $p:(S\setminus
O)\to P$. Let $\be\subset P$ be an infinite geodesic in direction
$\het$, and let $\al\subset S$ be its pull back by $p$. The
preceding argument shows that in the case i) (resp. case ii)) the
direction of $\al$ belongs to $U\setminus\E_{\erg}$ (resp.
$U\setminus(\E_{\erg}\cup\{\het_0\})$).
\end{proof}
\end{lem}

\medskip

Let $\tP$ be a noncompact, $G$-periodic polygonal surface, and let
$P=\tP/G$ be the compact quotient. Suppose that $\bo P\ne\emptyset$.
By Lemma~~\ref{bil_type_lem}, $U_{\bo P}P$ and $U_{\bo\tP}\tP$ are
the cross-sections for the billiard flows of $P$ and $\tP$
respectively. Let $(U_{\bo P}P,\tau,\nu)$ and
$(U_{\bo\tP}\tP,\ttau,\tnu)$ be the respective billiard maps. Then
$(U_{\bo\tP}\tP,\ttau,\tnu)$ is a skew product over $(U_{\bo
P}P,\tau,\nu)$ with the fibre $G$ and a displacement function $\vp:
U_{\bo P}P \to G$. See equation~~\eqref{produit_gauche_eq}.

\begin{lem}        \label{centr_lem}
The displacement function is centered:
\begin{equation}     \label{integr_eq}
\int_{U_{\bo P}P}\vp d\nu=0.
\end{equation}
\begin{proof}
For $v\in U_{\bo P}P$ let $\ga_v=\{\ga_v(t):0\le t\le t(v)\}$ be the
segment of the geodesic ray $\{\ga_v(t)\}$ determined by $v$ ending
when $\{\ga_v\}$ first returns to $\bo P$. Note that $ t(v)$ is the
{\em first return time function} for the billiard map. The tangent
vector $\ga'_v(t)$ is defined for $0<t<t(v)$ and the left limit
$\lim_{t\to t(v)-}\ga'_v(t)$ exists. Set $\si(v)=-\lim_{t\to
t(v)-}\ga'_v(t)$. The transformation $\si:U_{\bo P}P\to U_{\bo P}P$
is the {\em canonical involution} for the billiard map $(U_{\bo
P}P,\tau,\nu)$ \cite{Gut03}.

The canonical involution $\tsi$  for the billiard map
$(U_{\bo\tP}\tP,\ttau,\tnu)$ is defined the same way. Let
$\tv=(v,g)\in U_{\bo \tP}\tP$. Then
$$\tsi(\tv)=\tsi(v,g)=(\si(v),g+\vp(v)).$$
The identity
$$(v,g)=\tsi^2(v,g)=\tsi(\si(v),g+\vp(v))=(\si^2(v),g+\vp(v)+\vp(\si(v)))$$
yields
\begin{equation}     \label{symm_eq}
\vp(\si(v))=-\vp(v).
\end{equation}
Since $\si$ preserves the Liouville measure, the claim follows.
\end{proof}
\end{lem}

Let now $\tP$ be a rational polygonal surface. Let $\Ga(P)=R_N$; we
identify $U/\Ga(P)$ and $[0,\pi/N]$. For $\het\in[0,\pi/N]$ let
$(U_{\bo P}P_{\het},\tau_{\het},\mu_{\het})$ and
$(U_{\bo\tP}\tP_{\het},\ttau_{\het},\tmu_{\het})$ be the directional
billiard maps for $P$ and $\tP$ respectively. Then
$(U_{\bo\tP}\tP_{\het},\ttau_{\het},\tmu_{\het})$ is the skew
product over $(U_{\bo P}P_{\het},\tau_{\het},\mu_{\het})$ with the
displacement function $ \vp_{\het}=\vp|_{U_{\bo P}P_{\het}}$. By
equation~~\eqref{produit_gauche_eq}
\begin{equation}      \label{skew_prod_dir_eq}
\ttau_{\het}(v,g)=(\tau_{\het}(v),g+\vp_{\het}(v)).
\end{equation}

\begin{lem}      \label{integr_prop}
Let $N$ be even. Then for every $\het\in[0,\pi/N]$ the {\em
directional displacement function} $\vp_{\het}$ is centered:
\begin{equation}     \label{integr_eq1}
\int_{U_{\bo P}P_{\het}}\vp_{\het}d\nu_{\het}=0.
\end{equation}
If $N$ is odd, then the function $\vp_{\pi/(2N)}$ is centered.
\begin{proof}
The canonical involution $\si:U_{\bo P}P\to U_{\bo P}P$ induces {\em
directional involutions} $\si_{\het}:U_{\bo P}P_{\het}\to U_{\bo
P}P_{\eta(\het)}$. The central symmetry of $U$ and the
identification $U/R_N=[0,\pi/N]$ induce the transformation
$\het\mapsto\eta(\het)$ of $[0,\pi/N]$. The proof of
Lemma~~\ref{centr_lem} yields
\begin{equation}   \label{integr_eq2}
\int_{U_{\bo P}P_{\eta(\het)}}\vp_{\eta(\het)}d\nu_{\eta(\het)} =
-\int_{U_{\bo P}P_{\het}}\vp_{\het}d\nu_{\het}.
\end{equation}
If $N$ is even, then $R_N$ contains the central symmetry, hence
$\eta(\het)=\het$. If $N$ is odd, then $\eta(\het)=\pi/N-\het$. Both
claims now follow from equation~~\eqref{integr_eq2}.
\end{proof}
\end{lem}

\section{$\Z$-periodic polygonal surfaces}             \label{one_periodic}
We will use the setting and the notation of section~~\ref{cadre},
with $G=\Z$. Let $\tP$ be a $\Z$-periodic polygonal surface, and let
$P=\tP/\Z$. If $P$ is a rational polygonal surface, we will denote
by $\tS$ and $S$ the translation surfaces of $\tP$ and $P$
respectively. Then $\tS$ is $\Z$-periodic, and $S=\tS/\Z$.

\subsection{Main result}     \label{result_sub}

\hfill \break A compact translation surface $S$ is {\em arithmetic}
\cite{GJ} if it admits a translation covering $\pi:S\to\TT$ onto a
flat torus whose branch locus is a single point. Via an affine
renormalization, we can assume that $S$ covers the standard torus
$\TT_0=\RR/\ZZ$ and that the branch locus is $\{0\}+\ZZ$. These
translation surfaces are also known as {\em square-tiled} and as
{\em origamis}.


The surface $\tP$ is arithmetic iff $P=\tP/\Z$ is arithmetic
\cite{Gut09}. Let $P$ be a compact, arithmetic polygonal surface;
let $\pi:S\to\TT_0$ be as above. Let $\Ga=\Ga(P)$. A direction
$\het\in U$ is rational if $\tan\het\in\QQ$. Using the covering
$\pi:S\to\TT$ and the natural action of ${\text{GL}(2,\R)}$ on
translation surfaces \cite{GJ}, we extend the notion of {\em
rational directions} to all arithmetic translation surfaces, and
hence to arithmetic polygonal surfaces.\footnote{We will say {\em
$P$-rational} to emphasize that the set of rational directions
depends on the surface in question.}

Let $(U/\Ga)_{\rat}\subset U/\Ga$ be the set of $P$-rational
directions. Then $\het\in(U/\Ga)_{\rat}$ iff every geodesic  in $P$
in direction $\het$ is periodic or a {\em saddle connection}
\cite{Gut84}. The set $(U/\Ga)_{\rat}$ is countable. We set
$(U/\Ga)_{\irr}=U/\Ga\setminus(U/\Ga)_{\rat}$; we say that $\het\in
(U/\Ga)_{\irr}$ are the {\em irrational directions}.

\begin{thm}      \label{one_period_thm}
Let $\tP$ be a $\Z$-periodic polygonal surface with a boundary, and
let $P=\tP/\Z$.

\noindent 1. If the flow $(UP,T^t,\mu)$ is ergodic, then the
geodesic flow for $\tP$ is recurrent.

\noindent 2. Let $P$ be a rational polygonal surface, and let
$\Ga=\Ga(P)$. Suppose that $|\Ga|$ is divisible by $4$. Then for a
full measure set of directions $\het\in U/\Ga$ the directional
geodesic flow $(U\tP_{\het},\tT^t_{\het}, {\tla}_{\het})$ is
zero-recurrent. (See definition~~\ref{NRUO_def}.)

\noindent 3. Let $P$ be an arithmetic polygonal surface. i) For an
irrational direction $\het$ the flow $(U\tP_{\het},\tT^t_{\het},
{\tla}_{\het})$ is zero-recurrent. ii) Let $\het$ be a rational
direction. Then the set of orbits of $\tT_{\het}^t$ is a disjoint
union of periodic bands and bands of orbits that are dissipative
with a positive rate. The boundaries of these bands are
concatenations of saddle connections.
\begin{proof}
1. The claim follows from Lemma~~\ref{centr_lem} and claim 2 in
Proposition~~\ref{recurrence_prop}.

\noindent 2. Let $\het\in\E_{\erg}(P)$, the set of uniquely ergodic
directions. By Lemma~~\ref{bil_type_lem}, Lemma~~\ref{integr_prop},
and Proposition \ref{recurrence_prop}, the flow
$(U\tP_{\het},\tT^t_{\het}, {\tla}_{\het})$ is conservative. Since
the set $\E_{\erg}(P)\subset S^1$ has full lebesgue measure
\cite{KMS}, we obtain that $(U\tP_{\het},\tT^t_{\het},
{\tla}_{\het})$ is conservative for a.e. $\het$.

For $k\in\N$ set $P_k=\tP/k\Z$. Let $l,k\in\N$ and let $\ell$ divide
$k$. Then there is a covering $p_{k,\ell}:P_k\to P_\ell$, implying
$\E_{\erg}(P_k)\subseteq \E_{\erg}(P_\ell)$. In particular,
$\E_{\erg}(P_k)\subseteq \E_{\erg}(P)$ for any $k>1$. By
Lemma~~\ref{cobord_lem} and Proposition~~\ref{recurrence_prop}, if
$\het\in\E_{\erg}(P_k)$, and $k>1$, then the flow
$(U\tP_{\het},\tT^t_{\het}, {\tla}_{\het})$ is zero-recurrent. Thus,
$(U\tP_{\het},\tT^t_{\het}, {\tla}_{\het})$ is zero-recurrent for
$\het\in\cup_{k>1}\E_{\erg}(P_k)$, a full measure subset of
$\E_{\erg}(P)$.

\noindent 3. In this case all surfaces $P_k$ are arithmetic and
$\E_{\erg}(P_k)=\E_{\erg}(P)=(U/\Ga)_{\irr}$ \cite{Gut84}. We can
assume without loss of generality that
$\E_{\erg}(P)=[0,\pi/N]\setminus\QQ$. The preceding argument yields
claim i).


Let now $\het\in[0,\pi/N]\cap\QQ$. By \cite{Gut84}, the flow
$(UP,T_{\het}^t,\mu)$ decomposes into periodic bands whose
boundaries are made from saddle connections. Depending on whether
ergodic sums of the displacement function along a periodic orbit
vanish or not, the preimage of a periodic band in $U\tP$ is a union
of periodic and transient bands. Claim ii) follows.
\end{proof}
\end{thm}

\begin{corol}      \label{one_period_cor}
Let $\tP$ be a $\Z$-periodic, rational polygonal surface with a
boundary. If $|\Ga(P)|$ is divisible by $4$, then the flow
$(U\tP,\tT^t,\tmu)$ is zero-recurrent.
\begin{proof}
Follows from the decomposition of $(U\tP,\tT^t,\tmu)$ into the
directional flows $(U\tP_{\het},\tT^t_{\het}, {\tla}_{\het})$, a
Fubini-type argument, and claim 2 in Theorem~~\ref{one_period_thm}.
\end{proof}
\end{corol}

\subsection{Examples and applications}   \label{exa_sub}

\hfill \break  We will now illustrate the preceding material with a
few examples.

\begin{figure}[htbp]
\begin{center}
\input{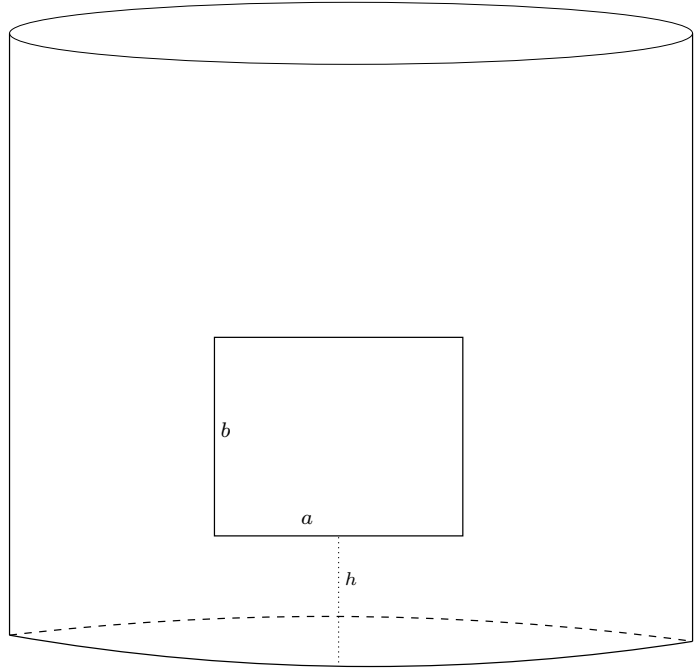}
\caption{\it Flat cylinder with a rectangular obstacle.}
\label{rect_obst}
\end{center}
\end{figure}

\begin{exa}      \label{more_strip_exa}
{\em Let $0<h<1$ and $0\le a,b <1$ be such that $0<h\pm\frac{b}2<1$
and $a+b>0$. Let $R(a,b;h)\subset\RR$ be the closed $a\times b$
rectangle centered at $(\frac12,h)$, whose sides are parallel to the
coordinate axes. Then $R(a,b;h)$ belongs to the interior of the unit
square $Q=\{(x,y):0\le x,y \le 1\}$. Let $P(a,b;h)$ be the polygonal
surface obtained by deleting from $Q$ the interior of $R(a,b;h)$,
and identifying the sides $\{x=0\},\{x=1\}$. If $0<a,b$, then
$P(a,b;h)$ is the flat unit cylinder with a rectangular obstacle.
The obstacle is the $a\times b$ rectangle centered in the cylinder
at the height $h$. See figure~~\ref{rect_obst}. If $b=0$ (resp.
$a=0$) then the rectangular obstacle degenerates into a horizontal
(resp. vertical) barrier.

For $k\in\Z$ let $R_k(a,b;h)=R(a,b;h)+(k,0)$; let
$B=\{(x,y):-\infty<x<\infty, 0 \le y \le 1\}$. Set
$\tP(a,b;h)=B\setminus\cup_{k\in\Z}R_k(a,b;h)$. Then $\tP(a,b;h)$ is
a $\Z$-periodic polygonal surface, and $P(a,b;h)=\tP(a,b;h)/\Z$.
When $h=\frac{1}2$, we recover Example~~\ref{strip_exa}.

Let $\Ga=\Ga(P(a,b;h))$. If $b\ne0$, then $|\Ga|=4$; when $b=0$,
then $|\Ga|=2$. Thus, for $b\ne 0$ the surface $\tP(a,b;h)$
satisfies the assumptions of claim 2 in
Theorem~~\ref{one_period_thm}. The surface $P(a,b;h)$ is arithmetic
iff $a,b\in\QQ$ \cite{GJ}. Theorem~~\ref{one_period_thm} and
Corollary~~\ref{one_period_cor} imply the following statement. }
\end{exa}

\begin{corol}       \label{rect_arithm_thm}
Let $(U\tP(a,b;h),\tT^t,\tmu)$ be the geodesic flow for
$\tP(a,b;h)$; let $(U\tP(a,b;h)_{\het},\tT^t_{\het},\tmu_{\het})$ be
the directional flows. We will refer to them as $\tT^t$ and
$\tT^t_{\het}$. Let $b\ne 0$. Then the following claims hold.

\noindent 1. The flow $\tT^t$ is zero-recurrent.

\noindent 2. For a.e. $\het\in[0,\pi/2]$ the flow $\tT^t_{\het}$ is
zero-recurrent.

\noindent 3. Let $a,b,h\in\QQ$. Then, for every $\het\in[0,\pi/2]$
such that $\tan\het\notin\QQ$, the flow $\tT^t_{\het}$ is
zero-recurrent.

\end{corol}

\begin{rem}       \label{horiz_bar_rem}
{\em If $b=0$, then $\tP(a,0;h)$ is the horizontal band with a
periodic configuration of horizontal barriers of length $a$. Thus
$N=1$ and $U/\Ga=[0,\pi]$. See Example~~15 in \cite{Gut09}. For
$\het\ne\pi/2$ the flows $\tT^t_{\het}$ are transient: Every orbit
of $\tT^t_{\het}$ drifts horizontally with the rate $\sin\het$. The
flow $\tT^t_{\pi/2}$ is periodic. This example fits into the
framework of Lemma~~\ref{integr_prop}. }
\end{rem}

\medskip

Let $Q\subset\RR$ be a polygon satisfying for an integer $t \geq 1$
the following conditions. i) There is a nonzero vector $\vv\in\RR$,
and for $1\le i \le t$ there are sides $s_i,s_i'$ of $Q$  such that
$s_i'=s_i+\vv$. ii) We have $Q\cap(Q+\vv)=\cup_{1\le i \le t}s_i'$.
We denote by $\tP=\tP(Q)$ the $\Z$-periodic polygon obtained by
deleting from $\cup_{k\in\Z}(Q+k\vv)$ the sides of the form
$s_i+k\vv:1\le i \le t,k\in\Z$. We say that $\tP$ is the {\em
stairway based on $Q$} or, simply, a {\em stairway}. The compact,
polygonal surface $P=\tP/\Z$ is obtained by identifying the sides
$s_i$ and $s_i'$ of $Q$ for $1\le i \le t$. Since $\cup_{1\le i \le
t}(s_i\cup s_i')\subset\bo Q$ is a proper subset, $\bo
P\ne\emptyset$. Let $\Ga\subset O(2)$ be the group generated by
reflections about the sides of $Q$ other than $s_i,s_i':1\le i \le
t$. If $|\Ga|<\infty$, then $\tP$ is a {\em rational stairway}. The
following is immediate from claims 2 and 3 in
Theorem~~\ref{one_period_thm}.

\begin{thm}        \label{stair_ration_thm}
Let $\tP\subset\RR$ be a rational stairway, and let $\Ga=R_N$. If
$N$ is even, then the following claims hold.

\noindent 1. For a full measure set of directions the flow
$(U\tP_{\het},\tT^t_{\het},\tmu_{\het})$ is zero-recurrent.

\noindent 2. Suppose, in addition, that the surface $\tP$ is
arithmetic. Then for every irrational direction the flow
$(U\tP_{\het},\tT^t_{\het},\tmu_{\het})$ is zero-recurrent.
\end{thm}

\begin{figure}[htbp]
\begin{center}
\input{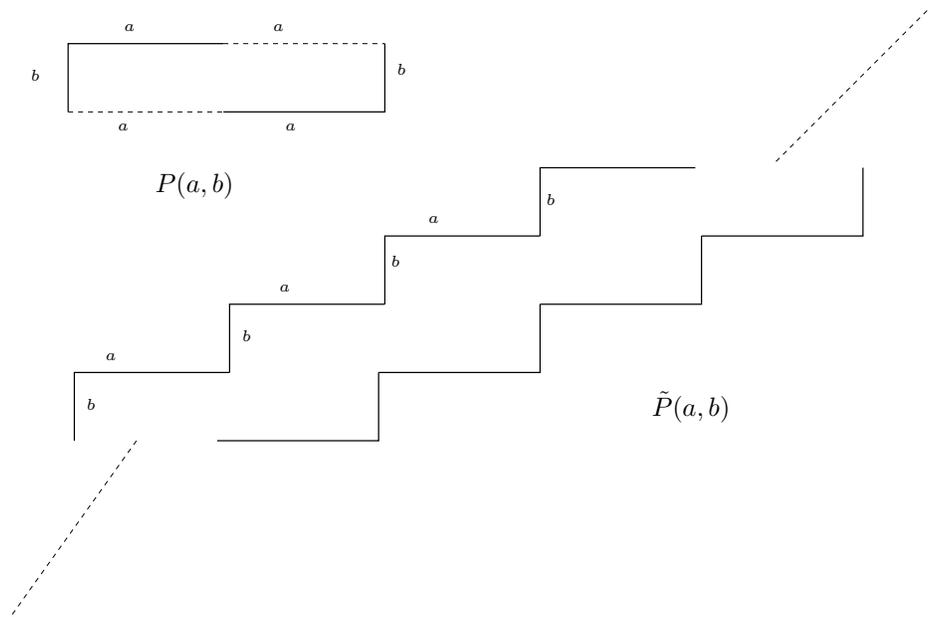}
\caption{\it A stairway polygonal surface and its quotient.}
\label{stairway}
\end{center}
\end{figure}

\begin{exa}     \label{stair_exa}
{\em Let $a,b>0$. Let $Q=Q(a,b)$ be the $2a\times b$ rectangle. We
view $\bo Q$ as a union of $6$ sides: 2 vertical sides of length $b$
and 4 horizontal sides  of length $a$. Let $s,s'$ be the lower left
and the upper right horizontal sides respectively. Then $\tP(a,b)$
based on $Q$ is the infinite stairway, with the stairs of length $a$
and height $b$. Its quotient $P(a,b)$ is the rectangle $Q$ with two
sides of length $a$ identified. See figure~~\ref{stairway}. The
corresponding group is $R_2$. By \cite{Gut84,GJ}, $\tP(a,b)$ is
arithmetic iff $a,b\in\QQ$.

}
\end{exa}

\begin{corol}          \label{stair_cor}
Let $\tP(a,b)$ be the stairway in Example~~\ref{stair_exa}. Then for
a.e. $\het\in[0,\pi/2]$ the flow $\tT^t_{\het}$ is zero-recurrent.
If $a,b\in\QQ$, then $\tT^t_{\het}$ is zero-recurrent if
$\tan\het\notin\QQ$.
\end{corol}

The claims of Corollary~~\ref{stair_cor} are immediate, by
Theorem~~\ref{stair_ration_thm}. See \cite{HoWe09} for another proof
of recurrence of $\tT^t_{\het}$ and \cite{HuWe09} for a study of
ergodic invariant measures.

\subsection{Generalizations and further applications}    \label{exten_sub}

\hfill \break Let $R_0$ be the unit square; let $O\subset R_0$ be a
polygon such that $R_0\setminus O$ is connected. Then
$\tP_O=\cup_{k\in\Z}\{(R_0\setminus\inter(O))+(k,0)\}$ is a periodic
band with obstacles and/or barriers. We will study the recurrence
for the billiard in $\tP_O$.

Let $R(a,b;h;\al)$ be the rectangle in Example~~\ref{more_strip_exa}
rotated by $\al$ about its center point. We assume that $a,b,h$ are
such that $R(a,b;h;\al)$ belongs to the interior of the unit square
for all $\al$. Set $\tP(\al)=\tP_{R(a,b;h;\al)}$.

Recall that a subset in a topological space is {\em residual} if it
contains a dense $G_{\de}$ set.

\begin{prop}        \label{rot_rect_thm}
The set of $\al\in S^1$ such that the flow $(U\tP(\al),\tT^t,\tmu)$
is recurrent is residual.
\begin{proof}
Set $P(\al)=\tP(\al)/\Z$. By \cite{KMS}, the billiard flow for
$P(\al)$ is ergodic for a dense $G_{\de}$ set of angles $\al$. The
statement now follows from claim 1 in Theorem~~\ref{one_period_thm}.
\end{proof}
\end{prop}

The set of planar polygons has a natural topology \cite{Gut03}. In
this topology, polygons with a fixed number of sides form closed
subsets in euclidean spaces. Imposing upper bounds on the sizes of
polygons, we obtain (relatively) compact subsets in euclidean
spaces. In what follows, whenever we invoke topological notions for
spaces of polygons, we mean the natural topology.

\begin{prop}        \label{tri_obst_thm}
Let $\ttt$ be the space of triangles in the interior of the unit
square. For $O\in\ttt$ let $\tP_O$ be the corresponding periodic
band with triangular obstacles. Then the set of triangles such that
the flow $(U\tP_{O},\tT^t,\tmu)$ is recurrent is residual.
\begin{proof} The space $\ttt$ is a relatively compact subset in $\R^6$.
For $O\in\ttt$ the quotient $P_O=\tP_O/\Z$ is the standard cylinder
with a triangular obstacle.  By \cite{KMS}, $\ttt$ contains a dense
$G_{\de}$ set of triangles such that the geodesic flow on $P_O$ is
ergodic. Now we apply claim 1 in Theorem~~\ref{one_period_thm}.
\end{proof}
\end{prop}

Let $\c$ be a closed set of polygons inside the unit square. For
$O\in\c$ let $\tP_O$ be the corresponding periodic band with
obstacles. Propositions~~\ref{rot_rect_thm} and~~\ref{tri_obst_thm}
are special cases of the following.

\begin{thm}     \label{dens_G_del_thm}
The set of $O\in\c$ such that the geodesic flow on $\tP_O$ is
recurrent is residual.
\end{thm}

The proof of Theorem~~\ref{dens_G_del_thm} is analogous to the
proofs of Propositions~~\ref{rot_rect_thm},~~\ref{tri_obst_thm}. It
invokes the result in \cite{KMS} that the ergodicity is
topologically typical. It is not known if ergodicity is typical
measure theoretically \cite{Gut03}.

\medskip

Y. Vorobets found a sufficient condition for the ergodicity of a
(compact) polygon \cite{Vo97}. The condition invokes the speed of
approximation of $\pi$-irrational angles of $P$ by rationals.
Referring the reader to \cite{Vo97} for a precise formulation, we
will say that the angles {\em admit the Vorobets approximation}.

\begin{thm}        \label{vorob_cond_thm}
Let $O$ be a polygon inside the unit square. Suppose that all
irrational angles of $O$ and those between $O$ and the horizontal
axis admit the Vorobets approximation. Suppose, moreover, that not
all of these angles are $\pi$-rational. Let $\tP_O$ be the
corresponding periodic band.  Then the geodesic flow on $\tP_O$ is
zero-recurrent.
\begin{proof}
For $k\in\N$ set $P_k=\tP_O/k\Z$. Since  all irrational angles of
$P_k$ admit the Vorobets approximation, the billiard flow on $P_k$
is ergodic. If $l$ divides $k$, we have the covering $p_{k,l}:P_k\to
P_l$. It remains to invoke the proof of claim 2 in
Theorem~~\ref{one_period_thm}.
\end{proof}
\end{thm}

We will now define a property of cocycles that has a simple
geometric meaning. Let $(X,\tau,\nu)$ be a dynamical system with a
finite invariant measure. Let  $\vp:X\to\R$ be a measurable
function, and let $(\vp_n)$ be the corresponding cocycle.

\begin{defin}     \label{unb_osc_def}
The cocycle $(\vp_n)$ has (the property of) {\it unbounded
oscillations} if for a. e. $x\in X$ we have
\begin{equation}   \label{unb_osc_eq}
\sup_n \vp(n,x) = +\infty, \ \inf_n \vp(n,x) = -\infty.
\end{equation}
\end{defin}

Let $(\tX,\tau_{\vp},\ttau)$ be the skew product over $(X,\tau,\nu)$
with the displacement function $\vp$. If the cocycle $(\varphi_n)$
has unbounded oscillations, then it is recurrent. Suppose that
$(\tX,\tau_{\vp},\ttau)$ is the Poincar\'e map of a skew product
flow. Then a. e. orbit of the flow has unbounded oscillations in the
obvious sense. See  \cite{BaKhMaPl} for examples of physical systems
corresponding to the billiard with unbounded oscillations  in
periodic polygons.

\begin{prop}    \label{unb_osc_thm}
Let $(X,\tau,\nu)$ be an ergodic dynamical system with finite
invariant measure. Let $\vp:X\to\R$ be a measurable function
satisfying $\int_X \vp \ d\nu = 0$; let $(\vp_n)$ be the
corresponding cocycle. If $\vp$ is not a coboundary, then $(\vp_n)$
has unbounded oscillations.
\begin{proof}
Suppose that the property $\sup_n\vp(n,x)=+\infty$ for a. e. $x\in
X$ is not satisfied. Then, by ergodicity, $\sup_n\vp(n,x)<\infty$
for a. e. $x\in X$. Set
\begin{equation}        \label{relat_eq}
h(x)= \sup_{k \ge 1} \vp(k, x), \  g(x) = \sup_{k \ge 2} \vp(k, x) -
h(x).
\end{equation}
Since $\tau \vp(k, x) = \vp(k+1, x) - \vp(x)$, we have
\begin{equation}   \label{phi_h_g_eq}
\vp(x)=\sup_{k \ge 2} \vp(k, x) - \tau \sup_{k \ge 1} \vp(k, x) =
h(x) - h(\tau x) + g(x).
\end{equation}
Iterating equation~~\eqref{phi_h_g_eq}, we obtain $\vp(n, x) = h(x)
- h(\tau^n x) + \sum_0^{n-1} g(\tau^jx)$.

\vskip 3mm By claim 2 in Proposition~~\ref{recurrence_prop}, the
cocycle $(\vp_n)$ is recurrent. Therefore,  for a.e. $x$ there is an
infinite sequence $n_k=n_k(x)$ such that $\vp (n_k, x)$ and
$h(\tau^{n_k} x)$ are bounded. Since, by equation~~\eqref{relat_eq},
$g(x)\le 0$, the above formula implies that the series
$\sum_{j=0}^\infty g(\tau^j x)$ converges for a.e. $x$. By the
recurrence of the cocycle, it implies $g=0$ a.e. Hence, by
equation~~\eqref{phi_h_g_eq}, $\vp$ is a coboundary, contrary to our
assumption. Assuming that the condition $\inf_n\vp_n(x)=-\infty$ for
a. e. $x\in X$ is not satisfied, we derive that $\vp$ is a
coboundary in a similar fashion.
\end{proof}
\end{prop}

Combining Proposition~~\ref{recurrence_prop} and
Proposition~~\ref{unb_osc_thm} with the statements on ergodicity in
\cite{KMS} and \cite{Vo97}, we strengthen the preceding results.
Below we formulate the strengthened versions of
Theorem~~\ref{dens_G_del_thm} and Theorem~~\ref{vorob_cond_thm}. We
leave the analogous strengthenings of
Proposition~~\ref{rot_rect_thm} and Proposition~~\ref{tri_obst_thm}
to the reader.

\begin{corol}      \label{dens_G_del_cor}
Let  $\c$ be a closed set of polygons inside  the unit square. For
$O\in\c$ let $\tP_O$ be the corresponding periodic band with
obstacles.

The set of $O\in\c$ such that the geodesic flow on $\tP_O$ is
zero-recurrent and has unbounded oscillations, is residual.
\end{corol}

\begin{corol}        \label{vorob_cond_cor}
Let $O$ be a polygon inside  the unit square; let $\tP_O$ be the
corresponding periodic band.

Suppose that each $\pi$-irrational angle of $O$ and each
$\pi$-irrational angle between $O$ and the horizontal axis admits
the Vorobets approximation. Suppose, moreover, that not all of these
angles are $\pi$-rational. Then the geodesic flow on $\tP_O$ is
zero-recurrent, with unbounded oscillations.
\end{corol}

We point out that there is a considerable interest in the physics
literature in the conservativeness and related properties for
$\Z$-periodic billiards. See, for instance, \cite{BaKhMaPl},
\cite{CLS}, and the references there.

\section{$\ZZ$-periodic polygonal surfaces: Dichotomies}    \label{two_periodic}

Let $R_0$ be the unit square; let $O\subset\inter(R_0)$ be a
polygon. Set $\tP_O=\RR\setminus\cup_{(p,q)\in\ZZ}(O+(p,q))$. Thus,
the $\ZZ$-periodic polygonal surface $\tP_O$ is the euclidean plane
with a doubly periodic configuration of obstacles. It is the
$\ZZ$-version of the band with obstacles studied in
section~~\ref{exa_sub} and section~~\ref{exten_sub}. The space
$\tP_O$ may be called the {\em polygonal Lorenz gas}. When $O$ is a
rectangle, $\tP_O$ is the wind-tree model. See  \cite{Eh} and
\cite{HW80}.

The study of geodesic flows for $\ZZ$-periodic polygonal surfaces is
less complete that the corresponding study for $\Z$-periodic
polygonal surfaces. In section~~\ref{rect_lorenz} we will study in
detail the directional flows for special directions in the wind-tree
model.

In this section we expose a few general results on the
conservativeness of arbitrary $\ZZ$-periodic polygonal surfaces
$\tP_O$.

For concreteness of exposition, we will consider the surfaces
$\tP_O$ when $O$ is a triangle. We may then call $\tP_O$ a {\em
triangular Lorenz gas}. The reader will easily extend the results
that follow to arbitrary polygons in the unit square. We denote by
$\ttt$ the topological space of triangles inside the unit square.

\begin{prop}                        \label{cas_irrat_thm}
There is a dense $G_{\de}$  set $\ddd\subset\ttt$ of triangles such
that for $O\in\ddd$ the following dichotomy holds: i) the billiard
in $\tP_{O}$ is zero-recurrent or ii) the billiard in $\tP_{O}$ is
transient and satisfies equation~~\eqref{lim-k_n}, with $d=2, c >
0$.
\begin{proof}
The billiard in $\tP_O$ fits into the framework of
sections~~\ref{trans_surf_sub},~~\ref{bill_sub}. The compact
polygonal surface $P_O=\tP_O/\ZZ$ is the standard torus with a
triangular obstacle. The flow $(U\tP_O,\tT^t,\tmu)$ is a skew
product over the geodesic flow $(UP_O,T^t,\mu)$ with the fibre
$\ZZ$. Set $\tO=\cup_{(p,q)\in\ZZ}(O+(p,q))$. The boundaries $\bo
O,\bo\tO$ yield canonical cross-sections for the respective flows.
Let $(U_{\bo O}P,\tau,\nu)$ and $(U_{\bo\tO}\tP,\ttau,\tnu)$ be  the
respective billiard maps. Then $(U_{\bo\tO}\tP,\ttau,\tnu)$ is a
skew product over $(U_{\bo O}P,\tau,\nu)$; let $\vp:U_{\bo
O}P\to\ZZ$ be the corresponding displacement function. The flow
$(U\tP_O,\tT^t,\tmu)$ is zero-recurrent (resp. transient) iff the
map $(U_{\bo\tO}\tP,\ttau,\tnu)$ is zero-recurrent (resp.
transient).

For $(p,q)\in\NN$ set $P_{(p,q)}=\tP_{O}/(p\Z\times q\Z)$. Then
$P_{(p,q)}$ are compact polygonal surfaces and $P_{(1,1)}=P_{O}$. If
$(p,q),(p',q')\in\NN$ are such that $p$ divides $p'$ and $q$ divides
$q'$, then there is a finite covering
$\pi_{(p,q)}^{(p',q')}:P_{(p',q')}\to P_{(p,q)}$. The set of
$O\in\ttt$ such that the billiard map $(U_{\bo
O}P_{(p,q)},\tau_{(p,q)},\nu_{(p,q)})$ is ergodic contains a dense
$G_{\de}$  set $\ddd_{(p,q)}$ \cite{KMS}.

For $O\in \cup_{(p>1,q)}\cup_{(p,q>1)}\ddd_{(p,q)}$, which is a
dense $G_{\de}$, the map $(U_{\bo\tO}\tP,\ttau,\tnu)$ satisfies the
conditions of Proposition~~\ref{lem-k_n}.
\end{proof}
\end{prop}

Let now $O\subset R_0$ be a rational triangle; thus, the polygonal
surface $P_O$ is rational; let $R_N$ be the corresponding dihedral
group. For $\het\in[0,\pi/N]$ let
$(UP_{\het},T^t_{\het},\mu_{\het})$ and
$(U\tP_{\het},\tT^t_{\het},\tmu_{\het})$ be the directional geodesic
flows on $P_O$ and $\tP_O$ respectively. Let $\tau_{\het}$ and
$\ttau_{\het}$ be the respective directional billiard maps. See
section~~\ref{bill_sub}.

\begin{prop}        \label{rat_case_thm}
Let $O\subset R_0$ be a rational triangle such that $N=N(O)$ is
even. Then for a.e. $\het\in[0,\pi/N]$ the following dichotomy
holds: i) The map $\ttau_{\het}$ is zero-recurrent or ii) the map
$\ttau_{\het}$ is transient and satisfies equation~~\eqref{lim-k_n},
with $d=2, c > 0$.
\begin{proof} Let $S,\tS$ be the translation surface
corresponding to $P_O,\tP_O$ respectively \cite{Gut84}. Then $\tS$
is a $\ZZ$-periodic translation surface, and $S=\tS/\ZZ$. We view
$T_{\het}^t,\tT_{\het}^t$ as flows on $S_O,\tS_O$ respectively.

For $(p,q)\in\N^2$ set $S_{(p,q)}=\tS/(p\Z\times q\Z)$. If
$(p,q),(p',q')\in\NN$ are such that $p$ divides $p'$ and $q$ divides
$q'$, then there is a finite covering
$g_{(p,q)}^{(p',q')}:S_{(p',q')}\to S_{(p,q)}$. We denote by
$(S_{(p,q)},T^t_{(p,q)},\mu_{(p,q)})$ the flow in direction $\het$
for the surface $S_{(p,q)}$. These flows are compatible with the
coverings $g_{(p,q)}^{(p',q')}$; they form a projective family.

Let $E\subset[0,\pi/N]$ be the set of directions $\het$ such that
all directional flows $T^t_{(p,q)}$ are ergodic. For $\het\in E$, by
Lemma~~\ref{bil_type_lem}, $\bo P_O$ yields cross-sections for the
flows $(S_{(p,q)},T^t_{(p,q)},\mu_{(p,q)})$. Let $\tau_{(p,q)}$ be
the Poincar\'e maps; let $\ttau$ be the Poincar\'e map with respect
to the corresponding cross-section for the flow
$(\tS,\tT^t_{\het},\tmu)$. Then $\ttau$ is a skew product over
$\tau_{(p,q)}$. By Lemma~~\ref{integr_prop}, the corresponding
displacement functions are centered. Proposition~~\ref{lem-k_n}
implies the claim for $\ttau$.

It remains to show that the set $E$ has full measure. For each
surface $S_{(p,q)}$ the set $E_{(p,q)}$ of uniquely ergodic
directions has full measure \cite{KMS}. Since
$E=\cap_{(p,q)\in\N^2}E_{(p,q)}$, the claim follows.
\end{proof}
\end{prop}

\section{$\ZZ$-periodic polygonal surfaces: Rectangular Lorenz gas}
\label{rect_lorenz}
We will study the polygonal Lorenz gas $\tP_O$ of
section~~\ref{two_periodic} when $O$ is a rectangle: The rectangular
Lorenz gas. In the physics literature this is known as the {\em
wind-tree model} \cite{HW80}; it is of some interest for foundations
of statistical physics \cite{Eh}. We begin by introducing notation.
Let $0< a, b< 1$. For $(m,n)\in\ZZ$ let $R_{(m,n)}(a,b)\subset\RR$
be the $a\times b$ rectangle centered at $(m,n)$ whose sides are
parallel to the coordinate axes. The Lorenz gas with rectangular
obstacles of size $a\times b$ corresponds to the polygonal surface
\begin{equation}   \label{lorenz_eq}
\tP(a,b) = \RR \setminus \bigcup_{(m,n)\in\ZZ}R_{(m,n)}(a,b).
\end{equation}
The quotient surface $P(a,b) = \tilde P(a,b)/\ZZ$ is the unit torus
with a rectangular hole.\footnote{The polygonal surface $\tP(0,b)$
or $\tP(a,0)$ is the euclidean plane with a $\ZZ$-periodic family of
barriers. The billiard flow is then transient. Hence, we assume that
$a,b> 0$.}

We will modify the notation of section~~\ref{cadre} as follows. The
surface $P(a,b)$ is rational and its dihedral group is $R_2$. We
identify $U/R_2$ with $[0,\pi/2]$. We will suppress $(a,b)$ from our
notation whenever this does not cause confusion. For
$\het\in[0,\pi/2]$ we denote by
$(\tZ_{\het},\tT^t_{\het},\tmu_{\het})$ and
$(Z_{\het},T^t_{\het},\mu_{\het})$ (resp.
$(\tX_{\het},\ttau_{\het},\tnu_{\het})$ and
$(X_{\het},\tau_{\het},\nu_{\het})$) the billiard flow (resp.
billiard map) for $\tP(a,b)$ and $P(a,b)$ respectively.

\subsection{Rational directions and small obstacles}
\label{small_obst_sub}
\hfill \break A direction $\theta\in[0,\pi/2]$ is rational if
$\tan\theta\in\QQ$. Rational directions $\theta(p,q) = arctan (q/p)$
correspond to pairs $(p,q)\in\NN$ with relatively prime $p,q$. When
there is no danger of confusion, we will use the notation  $(p,q)$
instead of $\het(p,q)$.

Let $R(a,b)=ABCD$ be the rectangle in the unit torus. Let
$\theta\in(0,\pi/2)$.  The space $X_{\het}$ consists of unit vectors
pointing outward, whose base points belong to $ABCD$ and whose
directions belong to the set $\{\pm\het,\pi\pm\het\}$. See
figure~~{\ref{bill-fig1}.

We say that the Lorenz gas $\tilde P(a,b)$ has {\it small obstacles
with respect to $(p,q)$} if the geodesics in $P(a,b)$ emanating from
$A$ or $C$ in the direction $\theta(p,q)$ return to either point
without encountering $R(a,b)$ on the way.

\begin{figure} [htbp]
\begin{center}
\includegraphics[scale=0.54]{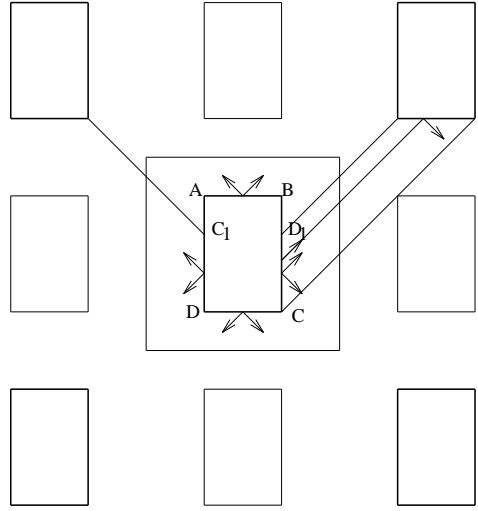}
\end{center}
\vspace{-2. cm} \caption{\it The cross-section for the conservative part of the billiard flow 
in direction $\pi/4$.}                      
\label{bill-fig1}
\end{figure}

\begin{lem} \label{small_obst_lem}
The small obstacles condition is satisfied iff
\begin{equation}
qa+pb \leq 1. \label{small_obst_eq}
\end{equation}
The inequality in equation~~\eqref{small_obst_eq} is strict iff the
directional geodesic flow $(\tZ_{(p,q)},\tT^t_{(p,q)},\tmu_{(p,q)})$
has a set of positive measure of orbits that do not encounter
obstacles.
\begin{proof}
The condition is satisfied iff $R(a,b)$ fits between two parallel
lines with slopes $q/p$ and vertical displacement $1/p$. By an
elementary calculation, this is possible iff
$$
a\frac{q}{p}+b\le\frac{1}{p}.
$$
Moreover, the equality in equation~~\eqref{small_obst_eq} holds iff
$R(a,b)$ takes all of the space between the boundary components of
the strip. See figure~~\ref{small_obst}.
\end{proof}
\end{lem}

\begin{figure}[htbp]
\begin{center}
\input{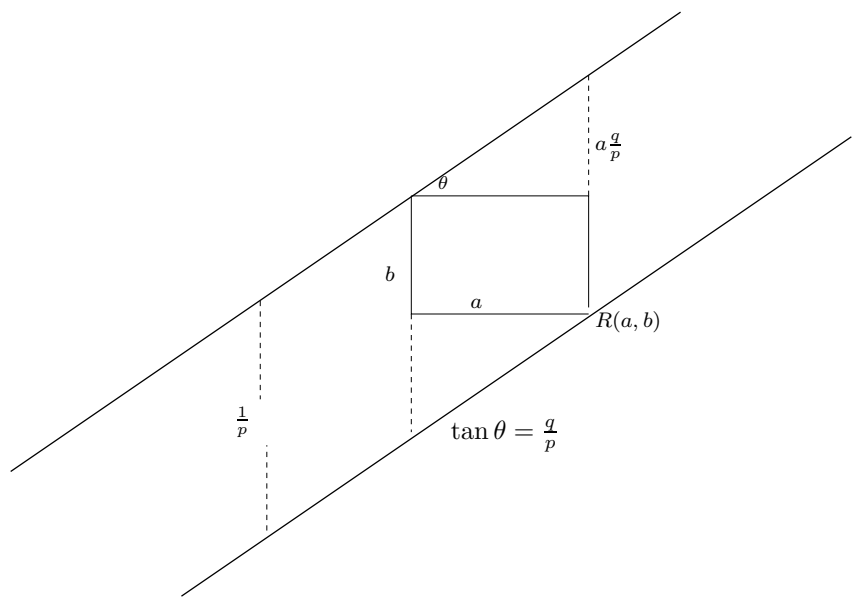}
\caption{\it The smallness of obstacles condition: Fitting the
rectangle in a strip.} \label{small_obst}
\end{center}
\end{figure}

In what follows we fix $(p,q)$ and assume that the inequality
equation~~\eqref{small_obst_eq} is satisfied. We identify
$X_{(p,q)}$ with 2 copies of the rectangle $ABCD$; the copy denoted
by $X_+=(ABCD)_+$ (resp. $X_-=(ABCD)_-$) carries the outward
pointing vectors in the directions $\het,\pi+\het$ (resp.
$\pi-\het,2\pi-\het$). Figure~~\ref{rectangles} illustrates this. We
will now investigate the Poincar\'e map $\tau_{(p,q)}:X_{(p,q)}\to
X_{(p,q)}$; we suppress the subscripts when it causes no confusion.

\begin{figure}[htbp]
\begin{center}
\input{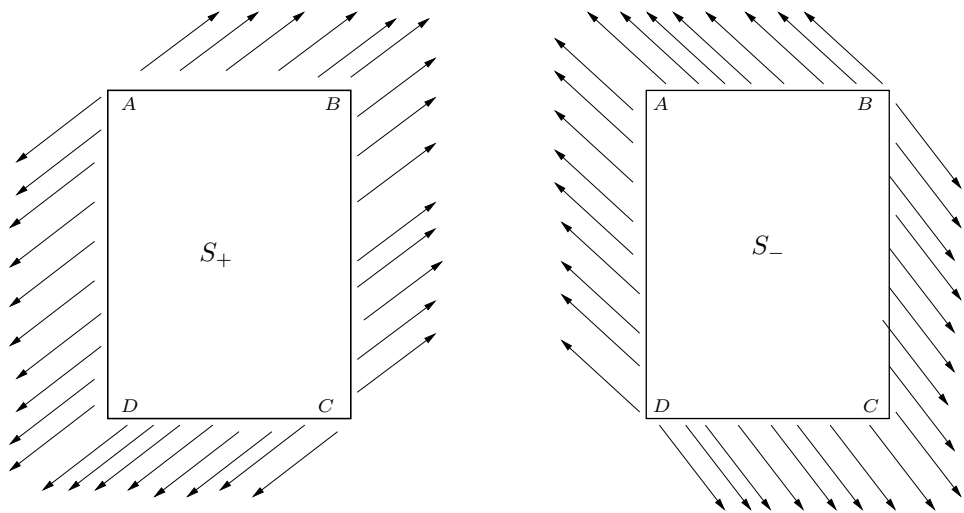}
\caption{\it The cross-section $X_{(p,q)}$ and the two corresponding
rectangles.} 
\label{rectangles}
\end{center}
\end{figure}

\begin{lem} \label{poinc_map_lem}
{\em 1.} There are natural identifications of  $X_+$ and $X_-$ with
the circle $\R/\Z$ endowed with distinguished points $A, B, C, D$;
their relative positions are given by
\begin{equation}   \label{arcs_eq1}
|AB| = |CD| = \frac{qa}{2(qa+pb)},\ |BC| = |DA| =
\frac{pb}{2(qa+pb)}.
\end{equation}
{\em 2.} Set $\tau_{\pm}=\tau|_{X_{\pm}}$. Then $\tau_+:X_+\to X_-,
\ \tau_-:X_-\to X_+$. Set $S=\R/\Z$. With the identifications
$X_{\pm}=S$, the maps $\tau_+:S\to S$ and $\tau_-:S\to S$ are the
orthogonal reflections about the axes $AC$ and $BD$ respectively.
The maps $\tau_-\tau_+:X_+\to X_+$ and $\tau_+\tau_-:X_-\to X_-$ are
the rotations of $S$ by $\displaystyle {qa\over (qa+pb)}$ and
$\displaystyle {pb\over (qa+pb)}$ respectively.
\begin{proof}
Vectors emanating from the rectangular obstacle in direction $\eta$
at the first return assume the direction $r(\eta)$, where $r$ is a
reflection in $R_2$. Thus, $\tau:X_+\to X_-,\ X_-\to X_+$.

Let $s$ and $\gamma$ be the arclength and the angle coordinates on
the billiard cross-section. Up to a constant factor, the invariant
measure for the billiard map has density $d\nu = \sin \gamma\, ds\,
d\gamma$.\footnote{See, e.g., \cite{GK} for this material.}
Integrating, we have $\nu(AB)=\nu(CD)=qa,\nu(BC)=\nu(DA)=pb$, up to
a constant factor. Normalizing $\nu(S)=1$, we obtain
equation~~\eqref{arcs_eq1}.

By construction, the maps $\tau_{\pm}:S\to S$ are orientation
reversing diffeomorphisms. Since they preserve the arclength, they
are isometries. Thus, $\tau_{\pm}:S\to S$ are orthogonal
reflections. By construction, $\tau_+$ (resp. $\tau_-$) fixes the
points $A,C$ (resp. $B,D$). These pairs of points correspond to the
axes of reflections when we identify $S$ with the unit circle in
$\RR$.
\end{proof}
\end{lem}

In what follows we will sometimes view $\tau_{\pm}$ as isometries of
the unit circle, $\tau_{\pm}:S\to S$, and sometimes as mappings
between the two copies of the circle, $\tau_+:S_+\to S_-, \
\tau_-:S_-\to S_+$.

Set $Z_{\pm}=S_{\pm}\times\ZZ$ and $Z=Z_+\cup Z_-$. We set
$\ttau=\ttau_{(p,q)}$. Then $\ttau:Z\to Z$ is the Poincar\'e map; it
interchanges the sets $Z_+,Z_-$. We use the notation
$\ttau_{\pm}(x,g)=(\tau_{\pm}(x),g+\vp_{\pm}(x))$. Thus,
$\vp_{\pm}:S\to\ZZ$ are the displacement functions. The following is
immediate from Lemma~~\ref{poinc_map_lem} and figure~~\ref{circle}.

\begin{figure}[htbp]
\begin{center}
\input{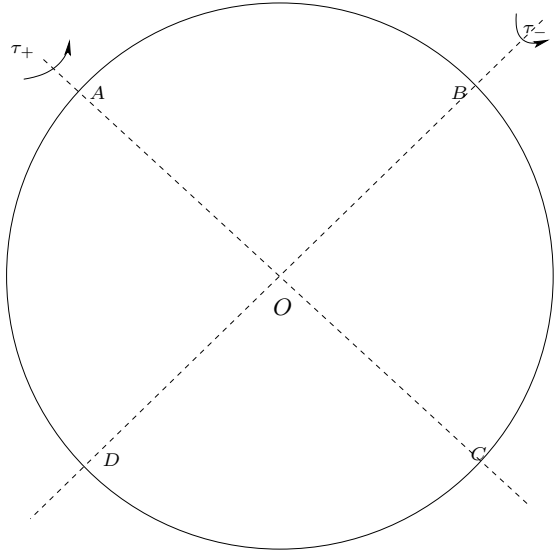}
\caption{\it The circle and the two orthogonal reflections.}
\label{circle}
\end{center}
\end{figure}

\begin{lem} \label{Poinc_map_lem}
The displacement functions $\vp_{\pm}:S\to\ZZ$ are constant on the
circular arcs $ABC$, $CDA$, $DAB$, $BCD$. We have
$$\vp_+|_{ABC}=(p,q),\ \vp_+|_{CDA}=(-p,-q),
\ \vp_-|_{DAB}=(-p,q),\ \vp_-|_{BCD}=(p,-q).$$
\end{lem}

Set
$$
A_1=\tau_-^{-1}(A),B_1=\tau_+^{-1}(B),C_1=\tau_-^{-1}(C),D_1=\tau_+^{-1}(D).
$$

Our next result describes the transformation $\ttau^2:Z\to Z$. We
set $\ttau^2_{\pm}=\ttau^2|_{Z_{\pm}}$. Recall that we have
identified $S$ and $\R/\Z$. We will usually denote by $x+y$ the
operation in $\R/\Z$. If the danger of confusion arises, we will
write $x+y\mod 1$.

\begin{prop} \label{poinc_map_prop}
We have
\begin{eqnarray*}
(\ttau^2)_+(x,g)&=&(x+\frac{qa}{qa+pb},g+\psi_+(x)),\\
(\ttau^2)_-(x,g)&=&(x+\frac{pb}{qa+pb},g+\psi_-(x)).
\end{eqnarray*}

The displacement functions $\psi_{\pm}$ take values $(\pm
2p,0),(0,\pm 2q)$. Each $\psi_{\pm}$ determines a partition of $S$
into four intervals such that $\psi_{\pm}=\const$ on each interval.
The endpoints of these intervals belong to the set
$A,B,C,D,A_1,B_1,C_1,D_1$.
\begin{proof}
We have $(\tau^2)_+=\tau_-\tau_+,(\tau^2)_-=\tau_+\tau_-$. The
product of two orthogonal reflections is the rotation by twice the
angle between their axes. The values of angles follow from
Lemma~~\ref{poinc_map_lem}.

We have $\psi_{\pm}(x)=\vp_{\pm}(x)+\vp_{\mp}(\tau_{\pm}(x))$. By
Lemma~~\ref{Poinc_map_lem}, $\psi_{\pm}$ are constant on the
circular arcs which are the intersections of half-circles
$ABC,CDA,DAB,BCD$ with half-circles
$A_1BC_1,C_1DA_1,D_1AB_1,B_1CD_1$. We have
$$
\psi_+|_{ABC\cap
D_1AB_1}=\vp_+|_{ABC}+\vp_-|_{DAB}=(p,q)+(-p,q)=(0,2q);
$$
$$
\psi_+|_{ABC\cap
B_1CD_1}=\vp_+|_{ABC}+\vp_-|_{BCD}=(p,q)+(p,-q)=(2p,0);
$$
$$
\psi_+|_{CDA\cap
D_1AB_1}=\vp_+|_{CDA}+\vp_-|_{DAB}=(-p,-q)+(-p,q)=(-2p,0);
$$
$$
\psi_+|_{CDA\cap
B_1CD_1}=\vp_+|_{CDA}+\vp_-|_{BCD}=(-p,-q)+(p,-q)=(0,-2q).
$$
Analogously
$$
\psi_-|_{DAB\cap
A_1BC_1}=\vp_-|_{DAB}+\vp_+|_{ABC}=(-p,q)+(p,q)=(0,2q);
$$
$$
\psi_-|_{DAB\cap
C_1DA_1}=\vp_-|_{DAB}+\vp_+|_{CDA}=(-p,q)+(-p,-q)=(-2p,0);
$$
$$
\psi_-|_{BCD\cap
A_1BC_1}=\vp_-|_{BCD}+\vp_+|_{ABC}=(p,-q)+(p,q)=(2p,0);
$$
$$
\psi_-|_{BCD\cap
C_1DA_1}=\vp_-|_{BCD}+\vp_+|_{CDA}=(p,-q)+(-p,-q)=(0,-2q).
$$
\end{proof}
\end{prop}

We set
\begin{equation}   \label{rot_angle_eq}
\alpha=\frac{qa}{qa+pb},\ \be=\frac{pb}{qa+pb}.
\end{equation}
Then $0<\alpha,\be<1$ and $\alpha + \be=1$. In what follows we
assume that $\al<\be$ or, equivalently, $qa<pb$. This assumption
allows us to avoid extra computations. The case $\be<\al$ reduces to
this by switching the coordinate axes. We  identify $S_+$ (resp.
$S_-$) with $[0,1]$ so that the points $A,B,C,D$ (resp. $D,A,B,C$)
go to $0,{\alpha\over 2},{1 \over 2},{1 \over 2}+{\alpha\over 2}$
(resp. $0,{1 \over 2}-{\alpha\over 2},{1 \over 2},1-{\alpha\over
2}$) respectively. With these identifications, $\psi_{\pm}:S\to\ZZ$
are piecewise constant functions on $[0,1]$. We will now explicitly
describe them. The formulas below follow from
Proposition~~\ref{poinc_map_prop} by straightforward calculations;
we leave them to the reader.
\begin{prop} \label{Poinc_map_prop}
The function $\psi_+:[0,1]\to\ZZ$ is given by
\begin{equation}   \label{cocy_plus_eq}
\psi_+(x) = \left \{ \begin{array}{clcr}
(0, 2q) &{\rm on \ }]0,{1\over 2}-{\alpha\over 2}[,\\
(2p, 0) & {\rm on \ } ]{1\over 2}-{\alpha\over 2},{1\over 2}[,\\
(0, -2q) & {\rm on \ } ]{1\over 2},1-{\alpha\over 2}[, \\
(-2p, 0) & {\rm on \ } ]1-{\alpha\over 2},1[.
\end{array} \right .
\end{equation}
The function $\psi_-:[0,1]\to\ZZ$ is given by
\begin{equation}   \label{cocy_minus_eq}
\psi_-(x) = \left \{ \begin{array}{clcr} (-2p, 0)
&{\rm on \ }]0,{\alpha\over 2}[,\\
(0, 2q) & {\rm on \ } ]{\alpha\over 2},{1\over 2}[,\\
(2p, 0) & {\rm on \ } ]{1\over 2},{1\over 2}+{\alpha\over 2}[, \\
(0, -2q) & {\rm on \ } ]{1\over 2}+{\alpha\over 2},1[.
\end{array} \right .
\end{equation}
\end{prop}

The following properties of $\psi_{\pm}:[0,1]\to\ZZ$ are immediate
from equations~~\eqref{cocy_plus_eq}, ~~\eqref{cocy_minus_eq}. The
lengths of intervals of continuity are ${\alpha \over 2},{\beta
\over 2}$, and they alternate. Each function takes four values which
generate the subgroup $H_{(p,q)}=2p\Z\oplus 2q\Z\subset\ZZ$. Using
the isomorphism $(a, b)\mapsto(2pa, 2qb)$ of $\Z^2$ and $H_{(p,q)}$,
we replace the displacement functions $\psi_+$ and $\psi_-$ by
piecewise constant functions on $[0,1]$ that do not depend on $p,q$.
Let $\Psi$ be the function corresponding to $\psi_+$. Then
\begin{equation}   \label{cocy_Psi}
\Psi(x) = \left \{ \begin{array}{clcr}
(0, 1) &{\rm on \ }]0,{1 \over 2}-{\alpha\over 2}[,\\
(1, 0) & {\rm on \ } ]{1 \over 2}-{\alpha\over 2},{1 \over 2}[,\\
(0, -1) & {\rm on \ } ]{1 \over 2},1-{\alpha\over 2}[, \\
(-1, 0) & {\rm on \ } ]1-{\alpha\over 2},1[.
\end{array} \right .
\end{equation}
\subsection{Ergodic decompositions for the billiard dynamics}   \label{decomp_sub}
\hfill \break Let $\ttau:X\times\ZZ\to X\times\ZZ$ be the billiard
map in direction $(p,q)$ for the Lorenz gas with rectangular
obstacles of size $a\times b$. Recall that we have identified $X$
with 2 copies of the unit circle: $X=S_+\cup S_-$. Let
$G_{(p,q)}\subset\ZZ$ be the group generated by $(p,q)$ and
$(p,-q)$. Then $|\Z^2/G_{(p,q)}|=2pq$. If $G$ is any countable
group, we will denote by $\tnu$ the measure on $X\times G$ which is
the product of the Lebesgue measure on $X$ and the counting measure
on $G$.

\begin{thm}     \label{erg_decom_thm}
Let $(p,q)\in\NN$ with $p,q$ relatively prime. Let $a,b>0$ satisfy
$qa+pb \leq 1$; suppose that $a/b$ is irrational. For
$\bg\in\Z^2/G_{(p,q)}$ denote by $\bg+G_{(p,q)}\subset\ZZ$ the
corresponding cosets. Then the following holds.

\noindent{\em 1.} For $\bg\in\Z^2/G_{(p,q)}$ the sets
$X\times(\bg+G_{(p,q)})\subset X\times\ZZ$ are $\ttau$-invariant.
The dynamical systems $(X\times(\bg+G_{(p,q)}),\ttau,\tnu)$ are
ergodic; they are isomorphic for all $\bg\in\Z^2/G_{(p,q)}$.

\noindent{\em 2.} The partition
\begin{equation}   \label{erg_dec_eq}
X\times\ZZ=\cup_{\bg\in\ZZ/G_{(p,q)}}X\times(\bg+G_{(p,q)}).
\end{equation}
yields the decomposition of the dynamical system
$(X\times\ZZ,\ttau,\tnu)$ into $2pq$ isomorphic ergodic components.
\begin{proof} Set $\bG=\ZZ/G_{(p,q)}$.
Lemma~~\ref{poinc_map_lem} and Lemma~~\ref{Poinc_map_lem} identify
$\ttau$ with the collection of transformations
$\ttau|_{\bg}:S_{\pm}\times(\bg+G_{(p,q)})\to
S_{\mp}\times(\bg+G_{(p,q)})$, where $\bg\in\bG$. This implies the
first part of claim 1.

Propositions~~\ref{poinc_map_prop} and~~\ref{Poinc_map_prop}
represent the restrictions
$\ttau^2_{\pm}|_{\bg}:S_{\pm}\times(\bg+G_{(p,q)})\to
S_{\pm}\times(\bg+G_{(p,q)})$ as skew product transformations
$\rho_{\si,\psi}$ over  certain rotations $s\mapsto s+\si$ on
$S=\R/\Z$ with particular displacement functions $\psi$. They do not
depend on $\bg\in\bG$. This proves the third part of claim 1.

With the notation of
equations~~\eqref{rot_angle_eq},~~\eqref{cocy_plus_eq},~~\eqref{cocy_minus_eq},
we have
\begin{equation}       \label{gbar_eq}
 \ttau^2_+|_{\bg}=\rho_{\al,\psi_+},\ \ttau^2_-|_{\bg}=\rho_{\be,\psi_-}.
\end{equation}

Since $\ttau$ interchanges $S_+\times(\bg+G_{(p,q)})$ and
$S_-\times(\bg+G_{(p,q)})$, the ergodicity of $\ttau|_{\bg}$ would
follow from the ergodicity of skew products
$\rho_{\al,\psi_+},\rho_{\be,\psi_-}$. By symmetry, it suffices to
prove the ergodicity of $\rho_{\al,\psi_+}$. Let $\Psi:S\to\ZZ$ be
given by equation~~\eqref{cocy_Psi}, and let
$\rho_{\al,\Psi}:S\times\ZZ\to S\times\ZZ$ be the corresponding skew
product. The isomorphism $G_{(p,q)}=\ZZ$ and
equation~~\eqref{cocy_Psi} yield
$\rho_{\al,\psi_+}=\rho_{\al,\Psi}$. By Theorem~~\ref{ergo-Psi} in
section~~\ref{spec_cocy_sub}, $\rho_{\al,\Psi}$ is ergodic for any
irrational $\al$. We have established claim 1. Claim 2 is immediate
from it.
\end{proof}
\end{thm}

Theorem~~\ref{erg_decom_thm} describes the ergodic decomposition of
the billiard map in direction $(p,q)$ on the polygonal surface
$\tP(a,b)$. We will now describe the decomposition of the geodesic
flow in direction $(p,q)$.\footnote{To simplify notation, we will
suppress the dependence on $(p,q)$ whenever this does not cause
confusion.} The configuration space for the directional flow
$(\tZ,\tT^t_{(p,q)},\tmu)$ consists of unit vectors in directions
$(\pm p,\pm q)$ with base points in $\tP(a,b)$. For $\tz\in\tZ$ we
denote by $\ga(\tz)\subset\tP(a,b)$ the geodesic it generates. Let
$\tC\subset\tZ$ (resp. $\tD\subset\tZ$) be the set of $\tz\in\tZ$
such that $\ga(\tz)$ encounters (resp. does not encounter)
rectangular obstacles. Then $\tZ=\tC\cup\tD$, a disjoint union. For
$\bg\in\bG$ set
$\tO(\bg)=\cup_{(m,n)\in(\bg+G_{(p,q)})}R_{(m,n)}(a,b)$. Thus,
$\tO(\bg)$ is the union of obstacles $R_{(m,n)}(a,b)$, as $(m,n)$
varies in the coset $\bg+G_{(p,q)}$. Let $\tC(\bg)\subset\tC$ be the
set of phase points $\tz\in\tZ$ such that $\ga(\tz)$ encounters
obstacles in $\tO(\bg)$. Let $\tmu_{\bg}$ be the restriction of
$\tmu$ to $\tC(\bg)$.

\begin{thm}    \label{bil_erg_decom_thm}
Let $(p,q)\in\NN$ with $p,q$ relatively prime. Let $a,b>0$
satisfying equation~~\eqref{small_obst_eq} be such that $a/b$ is
irrational. Let $(\tZ,\tT^t_{(p,q)},\tmu)$ be the directional flow.
Then the following holds.

\noindent{\em 1}. The sets $\tC(\bg),\bg\in\bG,$ are
$\tT^t$-invariant; the dynamical systems
$(\tC(\bg),\tT^t,\tmu_{\bg})$ are ergodic and pairwise isomorphic.
The partition $\tC=\cup_{\bg\in\bG}\tC(\bg)$ yields the ergodic
decomposition
\begin{equation}   \label{bil_erg_dec_eq}
(\tC,\tT^t_{(p,q)},\tmu)=\cup_{\bg\in\bG}(\tC(\bg),\tT^t_{(p,q)},\tmu_{\bg})
\end{equation}
of the conservative part of the flow $(\tZ,\tT^t_{(p,q)},\tmu)$.

\noindent{\em 2}.  The dissipative part $\tD$ is trivial iff we have
equality in equation~~\eqref{small_obst_eq}. Suppose that the
inequality in equation~~\eqref{small_obst_eq} holds. Then
$\tD=L\times\R$, where $L$ is a countable union of disjoint
intervals of the same length. The restriction of $\mu$ to $\tD$ is
the product of lebesgue measures on $L$ and $\R$; the flow
$(L\times\R,\tT^t_{(p,q)},\tmu)$ is the translation flow along $\R$.
\begin{proof}
By definition, the restriction of the flow $\tT^t_{(p,q)}$ to $\tD$
is dissipative. Claim 2 is immediate from
Lemma~~\ref{small_obst_lem}.

The flow $(\tC,\tT^t_{(p,q)},\tmu)$ is a suspension flow over the
transformation $(X\times\ZZ,\ttau,\tnu)$. Claim 1 now follows
directly from Theorem~~\ref{erg_decom_thm}. In particular,
equation~~\eqref{bil_erg_dec_eq} follows from the ergodic
decomposition of $(X\times\ZZ,\ttau,\tnu)$ given by
equation~~\eqref{erg_dec_eq}.
\end{proof}
\end{thm}

The proposition below relates the ergodic decomposition
equation~~\eqref{bil_erg_dec_eq} to an equidistribution of billiard
orbits. It holds under the assumptions of
Theorem~~\ref{bil_erg_decom_thm}.

The geodesic $\ga(\tz)$ generated by $\tz\in\tZ$ is a curve in the
polygonal surface $\tP(a,b)$. We will use the notation
$\ga_{\tz}(t),0\le t,$ for this curve, parameterized by the
arclength. For $(m,n)\in\ZZ,\tz\in\tC$ and $T>0$ let
$N(\tz,T;(m,n))$ be the number of times $0\le t \le T$ such that the
billiard orbit $\ga_{\tz}(t)$ encounters the obstacle
$R_{(m,n)}(a,b)$.

\begin{prop}    \label{bil_erg_decom_cor}
Let $(m,n),(m',n')\in\ZZ$. Then the following dichotomy holds.

\noindent 1. Suppose that the numbers
$\frac{m-m'}{p},\frac{n-n'}{q}$ are integers of the same parity.
Then there is a $\tT^t$ invariant subset $\tE\subset\tC$ of infinite
measure determined by the coset $(m,n)+G_{(p,q)}$, and such that for
$\tmu$-almost every $\tz\in\tE$ both functions
$N(\tz,T;(m,n)),N(\tz,T;(m',n'))$ go to infinity as $T\to\infty$.
Moreover, for $\tmu$-almost every $\tz\in\tE$ we have
\begin{equation}    \label{bil_erg_lim_eq}
\lim_{T \to \infty} \frac{N(\tz,T;(m,n))}{N(\tz,T;(m',n'))} = 1
\end{equation}
The set $\tC\setminus\tE$ also has infinite measure. For
$\tmu$-almost every $\tz\in\tC\setminus\tE$ we have
$$
N(\tz,T;(m,n))=N(\tz,T;(m',n'))=0.
$$

\noindent 2. Suppose that the above assumption on $(m,n),(m',n')$ is
not satisfied. Then for $\tmu$-almost every $\tz\in\tC$ one of the
following possibilities holds:

\medskip

\noindent a) $N(\tz,T;(m,n))=N(\tz,T;(m',n'))=0$;

\noindent b)  $N(\tz,T;(m,n))=0$, $N(\tz,T;(m',n'))\to\infty$;

\noindent c)  $N(\tz,T;(m',n'))=0$, $N(\tz,T;(m,n))\to\infty$.
\begin{proof}
1. Recall that $(\tC,\tT^t_{(p,q)},\tmu)$ is a suspension flow over
the billiard map $(X\times\ZZ,\ttau,\tnu)$. The set $X$ consists of
unit vectors with directions $(\pm p,\pm q)$ based on the boundary
of the rectangle $R(a,b)$. Let $\tX=X\times\ZZ$ and for
$(m,n)\in\ZZ$ set $\tX_{(m,n)}=X\times\{(m,n)\}$. Then
$$
\tX=\cup_{(m,n)\in\ZZ}\tX_{(m,n)},
$$
a disjoint union. We will use the ergodic decomposition of
$(X\times\ZZ,\ttau,\tnu)$ established in
Thorem~~\ref{erg_decom_thm}. For $\bg\in\bG$ let
$\tX(\bg)\subset\tX$ be the corresponding ergodic component of
$(X\times\ZZ,\ttau,\tnu)$. We denote by
$(m,n)\mapsto\overline{(m,n)}$ the projection of $\ZZ$ onto
$\bG=\ZZ/G_{(p,q)}$. Then, by equation~~\eqref{erg_dec_eq},
$$
\tX(\bg)=\cup_{\overline{(m,n)}=\bg}\tX_{(m,n)}.
$$
Denote by $1_{(m,n)}(\tx)$ the function $1_{\tX_{(m,n)}}:\tX\to\Z$.
For $k\in\N$ set
$$
f(\tx,k;(m,n))=\sum_{i=0}^k1_{(m,n)}(\ttau^i(\tx)).
$$

Our assumption on $(m,n),(m',n')$ is equivalent to
$\overline{(m,n)}=\overline{(m',n')}$. Set
$\overline{(m,n)}=\bg\in\bG$. Then
$\tX_{(m,n)},\tX_{(m',n')}\subset\tX(\bg)$. We use the ergodic
theorem for dynamical systems with infinite invariant measure
\cite{Aa97}. It states that for a.e. $\tx\in\tX(\bg)$
$$
\lim_{k \to \infty} \frac{f(\tx,k;(m,n))}{f(\tx,k;(m',n'))} =
\lim_{k \to \infty}
\frac{\sum_{i=0}^k1_{(m,n)}(\ttau^i(\tx))}{\sum_{i=0}^k1_{(m',n')}(\ttau^i(\tx))}=
\frac{\tnu(\tX_{(m,n)})}{\tnu(\tX_{(m',n')})}.
$$
The volume $\tnu(\tX_{(k,l)})$  does not depend on $(k,l)\in\ZZ$. We
have $0<\tnu(\tX_{(k,l)})<\infty$, and $\tnu(\tX_{(k,l)})$ is
determined by $a\times b$ and $(p,q)$. See
Lemma~~\ref{poinc_map_lem}. Hence, the preceding equation implies
the formula
\begin{equation}    \label{erg_lim_eq}
\lim_{k \to \infty} \frac{f(\tx,k;(m,n))}{f(\tx,k;(m',n'))} = 1,
\end{equation}
which holds, as usual, for a.e. $\tx\in\tX(\bg)$.

We will now outline an asymptotic relationship between the functions
$N(\tz,T;(m,n))$ and $f(\tx,k;(m,n))$. Let $P=\tP/\ZZ$ be the
compact polygonal surface; to simplify our notation, we suppress the
dependence on $a\times b$ and on $(p,q)$. Recall that $P$ is the
standard torus with a rectangular obstacle. Let $(Z,T^t,\mu)$ and
$(X,\tau,\nu)$ be the billiard flow and the billiard map in the
direction $(p,q)$ for $P$. Then $(Z,T^t,\mu)$ is a suspension flow
over $(X,\tau,\nu)$; the roof function $r(x):X\to\R_+$ is the time
it takes for the forward billiard orbit $\ga_x(t), 0<t,$ to return
to the cross-section. The mean time of return is given by
$$
\br=\frac{\int_Xr(x)d\nu}{\nu(X)}.
$$
Let $A\subset X$ be a measurable set. We associate with $A$ two
functions. The function $N_A(z,T):Z\times\R_+\to\N$ is the number of
times $0<t<T$ the billiard flow orbit $\ga_z(t)$ encounters $A$; the
function $f_A(x,k):X\times\N\to\N$ is the number of times $0\le i
<k$ the billiard map orbit $\tau^i(x)$ returns to $A$. Suppose that
$(X,\tau,\nu)$ is ergodic. Then for $\nu$-a. e. $x\in X$ we have
\begin{equation}    \label{flow_map_eq}
N_A(x,T) = f_A(x,\lfloor\frac{T}{\br}\rfloor)+o(T).
\end{equation}

We point out that equation~~\eqref{flow_map_eq} holds for ergodic
suspension flows, in general. In our setting the roof function has a
simple geometric meaning. Let $|\cdot|$ denote the euclidean norm on
$\RR$. Recall that $\vp:X\to\ZZ\subset\RR$ is the displacement
function. The elements $\tx\in\tX$ are vectors in $\RR$ based at
boundary points of the rectangles $R_{(m,n)}(a,b)$ as $(m,n)\in\ZZ$.
Let $\tx\mapsto x$ be the projection of $\tX$ onto $X$. For $x\in X$
let $b(x)\in\RR$ be the base point. Then $\tx\mapsto
b(\tau(x))-b(x)+\vp(x)$ is a well defined mapping of $\tX$ to $\RR$.
We then have
\begin{equation}    \label{bil_ret_eq}
r(\tx)=|b(\ttau(x))-b(x)+\vp(x)|.
\end{equation}
Note that the function $r:\tX\to\R_+$ is $\ZZ$-invariant. The mean
return time to the cross-section $\tX$ is equal to the mean return
time to the quotient cross-section $X$. We have
\begin{equation}    \label{mean_ret_eq}
\br=\frac{\int_X|b(\tau(x))-b(x)+\vp(x)|d\nu(x)}{\nu(X)}.
\end{equation}
Thus, $0<\br<\infty$. Set $\tE=\tC(\bg)$. Combining
equations~~\eqref{erg_lim_eq},~~\eqref{flow_map_eq},~~\eqref{mean_ret_eq},
we obtain the former part of our claim. We have
$\tC\setminus\tE=\cup_{\bh\in\bG\setminus\{\bg\}}\tC(\bh)$. The
remaining part of our claim follows from the preceding discussion
and equation~~\eqref{flow_map_eq}.

\medskip

\noindent 2. Set $\bg=\overline{(m,n)},\bg'=\overline{(m',n')}$. The
assumption on $(m,n),(m',n')$ does not hold iff $\bg\ne\bg'$. Thus,
$\tX(\bg)$ and $\tX(\bg')$ (resp. $\tC(\bg)$ and $\tC(\bg')$) are
distinct ergodic components of $\tX$ (resp. $\tC$). Equations a),
b), c) follow from the preceding discussion. Equation a) holds when
$\tz\in\tC\setminus(\tC(\bg)\cup\tC(\bg')$; equation b) (resp. c))
holds when $\tz\in\tC(\bg')$ (resp. $\tz\in\tC(\bg)$).
\end{proof}
\end{prop}

The following is immediate from
Proposition~~\ref{bil_erg_decom_cor}.

\begin{corol}                  \label{equidistr_cor}
Let $(m,n),(m',n')\in\ZZ$. Then for almost every $\tz\in\tC$ the
ratio $N(\tz,T;(m,n))/N(\tz,T;(m',n'))$ converges to either $1$, or
$0$, or infinity, as $T\to\infty$.
\end{corol}

\section{Ergodicity of cocycles over irrational rotations}  \label{ergodic}
The subject of this section is the ergodic theory for a certain
class of skew product transformations. The results are instrumental
in obtaining ergodic decompositions for directional flows in the
rectangular Lorenz gas model. See Theorem~~\ref{erg_decom_thm},
Theorem~~\ref{bil_erg_decom_thm}, and
Corollary~\ref{bil_erg_decom_cor}.

Throughout this section, we will use the following setting. Let $G$
be a locally compact abelian group.\footnote{In our applications
$G\simeq\R^{d_1}\times\Z^{d_2}$.} Set $X = \R/\Z$. For $0<\al<1$ let
$\rho_\alpha: X\to X$ be the rotation $x \mapsto x+ \alpha {\rm \
mod \ 1 \ }$. Let $\Phi:X\to G$ be a piecewise constant function.
Define the transformation $\rho_{\alpha,\Phi}:X\times G\to X\times
G$ by $(x,g) \mapsto (\rho_\alpha(x), g+\Phi(x))$. Let $\leb$ denote
the Lebesgue measure on $X$; let $\mu$ be the measure on $X\times G$
which is the cartesian product of $\leb$ and a Haar measure on $G$.
The dynamical system $(X\times G,\rho_{\alpha,\Phi},\mu)$ is the
{\em skew product} over $\rho_\alpha$, with the fibre $G$ and the
displacement function $\Phi$. Let $(\Phi_n)$ be the cocycle
corresponding to $\Phi$ and $\rho_{\al}$. See
Definition~~\ref{cocycle_def}. If the dynamical system $(X\times
G,\rho_{\alpha,\Phi},\mu)$ is ergodic, we will say that {\em the
cocycle $(\Phi_n)$ is ergodic}.

In the studies of ergodicity for $(X\times G, \rho_{\alpha, \Phi},
\mu)$ it is common to assume that the points of discontinuity of
$\Phi$ are not arithmetically related to $\alpha$. In view of
equation~~\eqref{cocy_Psi}, this assumption does not hold for our
applications in section~~\ref{rect_lorenz}. Thus, we will expose two
approaches to proving the ergodicity  of cocycles over irrational
rotations. One of them is based on the ``well distributed
discontinuities'' property for a cocycle.\footnote{We will use the
abbreviation (wdd).} See Definition~~\ref{propert_def}. The cocycles
$(\Psi_n)$, needed for our applications, satisfy (wdd) for generic
$\al$. This approach allows us to establish the ergodicity of
$(\Psi_n)$, and similar cocycles, for generic rotation angles. We
develop this approach in section~~\ref{erg_cocy_sub}, see especially
Corollary~~\ref{D-wdd}.

In our applications in section~~\ref{rect_lorenz}, $\al$ is
determined by the parameters $a$ and $b$, i. e., the sizes of
billiard obstacles. See equation~~\eqref{rot_angle_eq}. Hence, the
results of section~~\ref{erg_cocy_sub} prove the claims of
section~~\ref{decomp_sub} for generic obstacles. However, (wdd) may
fail for some parameters $a,b$. Our other approach is geared
specifically to the cocycles $(\Psi_n)$. We establish their
ergodicity for all irrational $\al$ in section~~\ref{spec_cocy_sub}.

\subsection{The inequality of Denjoy-Koksma}        \label{denj_koks_sub}
\hfill \break We recall basic facts about continued fractions. See,
for instance, \cite{Kh37} for this material. Let $\alpha \in ]0,1[$
be an irrational number, let $[0;a_1,..., a_n,...]$ be its continued
fraction representation, and let $(p_n/q_n)_{n \ge 0}$ be the
sequence of its convergents. The integers $p_n$ (resp. $q_n$) are
the {\em numerators} (resp. {\em denominators}) of $\alpha$. Thus
$p_{-1}=1$, $p_0=0$, $q_{-1}=0$, $q_0=1$. For $n \ge 1$ we have
\begin{equation} \label{converg_eq}
p_n = a_n p_{n-1}+p_{n-2}, q_n = a_n q_{n-1}+q_{n-2},  (-1)^n =
p_{n-1} q_n - p_n q_{n-1}.
\end{equation}
For $u\in\R$ set $\|u\|= \inf_{n \in \Z} |u - n|$. Then for $n \ge
0$ we have $\|q_n \alpha\| = (-1)^n (q_n \alpha - p_n)$. We also
have
\begin{eqnarray}
1 &=& q_n\|q_{n+1} \alpha \| + q_{n+1} \|q_n \alpha\|, \label{f_1} \\
{1\over q_{n+1}+q_n} &\le& \|q_n \alpha\| \le {1\over q_{n+1}}
= {1 \over a_{n+1} q_n+q_{n-1}}, \label{f_3} \\
\|q_n \alpha \| &\leq& \|k \alpha \|\ \mbox{for}\  1\le k < q_{n+1}
\label{f_4}.
\end{eqnarray}

We denote by $V(\varphi)$ the variation of functions; we will use
the shorthand BV for functions of  bounded variation. A function is
{\em centered} if $\int_X\vp(x)\,dx = 0$. Let $\varphi$ be a
centered BV function on $X$. Let $p/q$ be a rational number in
lowest terms such that $\|\alpha - p/q\| < {1 / q^2}$. The {\em
Denjoy-Koksma inequality} says that for any $x\in X$ we have
\begin{eqnarray}
|\sum_{\ell = 0}^{q-1}\varphi(x+\ell \alpha)| \le V(\varphi).
\label{f_8}
\end{eqnarray}

The following is immediate from equations~~\eqref{f_3} and
\eqref{f_8}.

\begin{cor}     \label{rec1}
Let $\Phi:X\to G$ be a centered BV function. Then the cocycle
$(\Phi_n)$ over any irrational rotation is recurrent.
\end{cor}

We will use the following properties of the sequence $p_k/q_k$. At
least one of any two consecutive numerators (resp. denominators) is
odd. If both $p_n, q_n$ are odd, then one of $p_{n+1}, q_{n+1}$ is
even.

\begin{lem}     \label{4odd}
Let $\al$ and $p_k,q_k$ for $k\ge 0$ be as above. Then the following
holds. {\em 1}. In any pair of  consecutive denominators at least
one satisfies $q_n\|q_n\alpha\| < 1/2$. {\em 2}. Out of any four
consecutive denominators at least one is odd and satisfies $q_n
\|q_n \alpha\| < 1/2$.
\begin{proof}
1. For any $n\in\N$ define $\delta_1$ and $\delta_2$ by
$q_n\|q_n\alpha\| = {1\over 2} - \delta_1$,
$q_{n+1}\|q_{n+1}\alpha\| = {1\over 2} - \delta_2$. By
equation~~\eqref{f_1}, we have
$$(q_{n+1} - q_n)^2 = 2\delta_1 q_{n+1}^2 + 2\delta_2q_{n}^2.$$
Hence, $\delta_1$ and $\delta_2$ cannot be both negative.

\noindent 2. The a priori possible parities for any four consecutive
denominators $q_{n-1}, q_n, q_{n+1}, q_{n+2}$ are as follows:
\begin{eqnarray*}
&&(0, \ 1, \  0, \  1) \ \ (0, \  1, \  1, \  0)  \ \ (1, \  0, \
1, \  0) \ \ (0, \  1, \ 1,
\  1)  \\
&&(1, \  0, \  1, \  1) \ \ (1, \  1, \  0, \  1) \ \ (1, \  1, \ 1,
\  0) \ \ (1, \  1, \  1, \  1).
\end{eqnarray*}
If there are two consecutive odd denominators, then the statement
follows from claim 1. It remains to consider the possibilities $(0,
\ 1, \ 0, \ 1)$ and $(1, \ 0, \  1, \  0)$. Then we have,
respectively, $q_n$ is odd, $a_{n+1} \not =1$, and $q_{n+1}$ is odd,
$a_{n+2} \not =1$. Set $q=q_{n}$ (resp. $q = q_{n+1}$) in the former
(resp. latter) case. Then $q\|q\alpha\| < 1/2$.
\end{proof}
\end{lem}

\subsection{Ergodicity of generic cocycles} \label{erg_cocy_sub}
\hfill \break Let $d(\cdot,\cdot)$ be an invariant distance on $G$.
\begin{defin}            \label{def-period}
{\em Let $a\in G$. 1. Suppose that for $n \ge 1$  there exist
$\ell_n\in\N$, $\varepsilon_n > 0,$ and $\delta>0$ such that
\begin{itemize}
\item[i)] We have $\lim_n \varepsilon_n = 0$, $\lim_n \ell_n \alpha {\rm \ mod \ 1 \
} = 0$,
\item[ii)] We have $\leb(\{ x: d(\Phi_{\ell_n}(x) , a) < \varepsilon_n \}) \ge \delta$.
\end{itemize}
Then we say that $a$ is a {\it quasi-period} for the cocycle
$(\Phi_n)$.

\noindent 2. We say that $a$ is a {\it period} if for every
$\rho_{\alpha, \Phi}$-invariant measurable function $f$ on $X\times
G$ and for a.e. $(x,g) \in X \times G$ we have
\begin{equation}  \label{period_eq}
f(x,g+a) = f(x,g).
\end{equation}

}
\end{defin}

We will use the following fact  \cite{Co09}.

\begin{lem}  \label{lem-period}
Every quasi-period is a period.
\end{lem}

The set of periods is a closed subgroup of $G$ which coincides with
the group of finite essential values of the cocycle \cite{Sc77}. A
cocycle is ergodic iff its group of periods is $G$.

\medskip

We introduce more notation. Let $\Phi:X\to G$ be a non constant
piecewise constant function. Denote by $R(\Phi)\subset G$ the range
of $\Phi$, i. e., $a\in R(\Phi)$ iff $\Phi(x) =a$ on a nontrivial
interval. Denote by ${\mathcal D} = \{t_i: i=1,...,d\}$ the set of
discontinuities of $\Phi$. We assume without loss of generality that
$0\in{\mathcal D}$. For $N \in \N$ let ${\mathcal D}_N=\{t_i - j\al
{\rm \mod 1}:\,1\le i\le d,0\le j < N\}$ be the set of
discontinuities for $\Phi_N(t)= \sum_{k=0}^{N-1} \Phi(t + k\alpha)$.
We set ${\mathcal D}_N=\{0= \gamma_{N,1} < ... <\gamma_{N,dN}<1\}$;
thus, for $1\le\ell\le dN$ the elements $\ga_{N,\ell}$ run through
${\mathcal D}_N$ in the natural order. We set $\gamma_{N,dN+1} =
\gamma_{N,1}$. The following notions will be important in what
follows.

\begin{defin}    \label{propert_def}
{\em 1. Let $0< \alpha <1$ be irrational. Let $\Phi:X\to G$ be a
piecewise constant function; let $(\Phi_n)$ be the corresponding
cocycle over $\rho_\alpha$. Suppose that there is $c >0$ and an
infinite set $W$ of denominators of $\alpha$ such that for all $q
\in W,\ell \in \{1,\dots,dq\}$ we have
\begin{eqnarray}  \label{wddPrty}
\gamma_{q,\ell+1} - \gamma_{q,\ell} \geq {c \over q}.
\end{eqnarray}
Then the cocycle has {\em well distributed discontinuities}. We
will use the shorthand {\em (wdd)}.

\noindent 2. Let $\alpha \in ]0,1[$ be irrational, and let
$[0;a_1,..., a_n,...]$ be its continued fraction. Then  $\al$ has
{\em property (D)} if there is $M\in\N$ such that for \ infinitely
many $n$ either $a_{n} \in [2, M]$ or $a_{n}= a_{n+1} =1$. }
\end{defin}

\begin{lem} \label{min-ri2}
Let $t\in\frac12(\Z\al+\Z) \setminus (\Z\al+\Z)$. Then there exist
$c>0$ and $L\in\N$ such that, if $n \geq L$ and either i) $a_{n+1}
\in [2, M]$ or ii) $a_{n+1} = a_{n+2} = 1$, then for $1\le k \le
q_n-1$ we have

\begin{equation}   \label{twell-dist}
\|k \alpha - t\| \geq {c\over q_n}.
\end{equation}
\begin{proof}
Let $t= {1\over 2} \ell \alpha + {1\over 2}r$, with $\ell, r \in
\Z$, and $\ell$ or $r$ odd. Let $L$ be such that $|\ell|  < q_{n-1}$
for $n \geq L$. Let $k \in [1, q_n-1]$. For $n \geq L$, we have if
$a_{n+1} \geq 2$,
$$ |2k - \ell| \leq  2k + |\ell| < 2 q_n+q_{n-1} \leq a_{n+1}q_n+q_{n-1} = q_{n+1}.$$

If $a_{n+1} \in [2, M]$, then, by equations~~\eqref{f_3}
and~~\eqref{f_4}, for all $j \in [1, q_{n+1}[$, we have
\begin{equation}    \label{minorJ}
\|j \alpha \| \geq \|q_n \alpha \| \geq {1\over q_n+q_{n+1}} =
{1\over (a_{n+1}+1)q_n+q_{n-1}} \geq {1\over (2+M)q_n}.
\end{equation}

If $2k - \ell \not = 0$, equation~~\eqref{minorJ} implies
$$\|k \alpha - t\| \geq {1\over 2} \|(2k - \ell)\alpha\|
\geq \frac{1}{2q_n(2 + M)}.$$

If $a_{n+1} = a_{n+2} = 1$, we have $q_{n+1} = q_{n}+ q_{n-1}$,
$q_{n+2} = q_{n+1}+ q_n = 2 q_n +q_{n-1}$ and $|2k - \ell| < 2
q_n+q_{n-1} = q_{n+2}$. Thus, if $2k - \ell \not = 0$
$$\|(2k - \ell) \alpha\| \geq \|q_{n+1} \alpha\| \geq
{1\over q_{n+2}+q_{n+1}} = {1\over 3q_{n}+2q_{n-1}} \geq {1 \over 5
q_{n}}.$$

\vskip 3mm If $\ell$ is even then $r$ is odd, so that for $2k =
\ell$ we have $\|k \alpha - t\| = \|(k-{1\over 2}\ell) \alpha -
{1\over 2}\| = {1\over 2}$.

\vskip 3mm In each case equation~~\eqref{twell-dist} holds with $c=
\inf(\frac{1}{2(2 + M)}, \frac1{10})$.
\end{proof}
\end{lem}

\vskip 3mm The following is immediate from Lemma~~\ref{min-ri2}.

\begin{prop} \label{wdd-ri2}
Let $\alpha$ be irrational, let $\Phi:X\to G$ be a piecewise
constant function, and let ${\mathcal D}\subset X$ be its set of
discontinuities. If $\alpha$ has property (D) and if ${\mathcal
D}\subset\{\frac12(\Z\al+\Z)  \setminus (\Z\al+\Z)\}\mod 1$, then
the cocycle $(\Phi_n)$ has (wdd).
\end{prop}

The following theorem is the main result of this subsection.

\begin{thm}          \label{wdd-erg-thm}
Let $0<\al<1$ be irrational. Let $\Phi:X\to\Z^r$ be a piecewise
constant, centered function, and let $R(\Phi)\subset\Z^r$ be the
range of $\Phi$.

Let $(\Phi_n)$ be the corresponding cocycle. If $(\Phi_n)$ has
property (wdd), then the group of its periods contains $R(\Phi)$.
\begin{proof} 
Let $W\subset\N$ be the infinite set of denominators of $\al$
introduced in Definition~~\ref{propert_def}. We will study the
family of functions $\{\Phi_q(x)=\sum_{0\le k \le
q-1}\Phi(x+k\al\mod 1)\,:q\in W\}$. For $1\le \ell \le dq$ let
$I_{q,\ell}= ]\gamma_{q,\ell}, \gamma_{q,\ell+1}[$ be the intervals
of continuity. Set $\II_q=\{I_{q,\ell}:1\le \ell \le dq\}$.

For $t_i$ in the set ${\mathcal D}$
of discontinuities of $\Phi$, let $\sigma_i=\lim_{\varepsilon
\rightarrow 0^+} [\Phi(t_i +\varepsilon) - \Phi(t_i - \varepsilon)]$
be the jump at $t_i$; let $\Si(\Phi)=\{\si_i:1\le i \le d\}$ be the
set of jumps. Set $R=\cup_{q\in W}R(\Phi_q) \subset \Z^r$. By
equation~~\eqref{f_8}, $R$ is a finite set.

\medskip

Each interval $[{k \over q}, {k+1 \over q}[$, $0 \leq k < q$,
contains an element $j\alpha \mod 1, 0\le j \le q-1$. Thus, for any
$t\in X$, the elements $\{t+j\alpha\mod 1:j=0,\dots,q-1\}$ partition
$X$ into intervals of lengths less than $2/q$. Hence, any interval
$J\subset X$ of length $\ge 2/q$ contains at least one point of the
set $\{t+j\alpha\mod 1:j=0,\dots,q-1\}$.

Let $c$ be the constant from equation~~\eqref{wddPrty}; set $L =
\lfloor\frac{2}{c}\rfloor+1$. Let $I_{q,\ell} \in \II_q$ be
arbitrary. Set $J_{q, \ell} \subset X$ be the union of $L$
consecutive intervals in $\II_q$ starting with $I_{q,\ell}$. By
equation~~\eqref{wddPrty}, the length of $J_{q,\ell}$ is greater
than or equal to $2/q$. Thus, for any $t_i\in\ddd$, the interval
$J_{q, \ell}$ contains a point in the set $\{t_i+j\alpha\mod
1:j=0,\dots,q-1\}$. Therefore, for any $\sigma \in \Sigma(\Phi)$,
there is $v\in R$ and two consecutive intervals $I,I'\in\II_q$ such
that $I\cup I'\subset J_{q, \ell}$ and such that $\Phi_q$ takes
values $v$ and $v+\si$ on $I$ and $I'$ respectively.

Let $\sigma \in \Sigma(\Phi)$, $v\in R$. Let
$\FF_q(\sigma)\subset\II_q$ be the family of intervals $I\in\II_q$
such that the jump at the right endpoint of $I$ is $\sigma$. Let
$\A_q(\sigma,v)\subset\FF_q(\sigma)$ be the set of intervals
$I\in\FF_q(\sigma)$ such that the value of $\Phi_q$ on $I$ is $v$;
let $\A'_q(\sigma,v)\subset\FF_q(\sigma)$ be the set of intervals
$I'\in\FF_q(\sigma)$ adjacent on the right to the intervals
$I\in\A_q(\sigma,v)$. Let $A_q(\sigma,v)\subset X$ (resp.
$A'_q(\sigma,v)\subset X$) be the union of intervals
$I\in\A_q(\sigma,v)$ (resp. $I'\in\A'_q(\sigma,v)$). Thus, $\Phi_q$
takes value $v$ (resp. $v+\si$) on $A_q(\sigma,v)$ (resp.
$A'_q(\sigma,v)$).

Denote by $|\cdot|$ the cardinality of a set. There is $v_0 \in R$
and an infinite subset of $W$ (which we denote by $W$ again) such
that for $q\in W$ we have
\begin{equation}   \label{proport_eq}
|\A_q(\sigma,v_0)|,|\A'_q(\sigma,v_0)|\ge\frac{1}{|R|}|\FF_q(\sigma)|.
\end{equation}
We have $|\FF_q(\sigma)|\ge qd/L$. By equation~~\eqref{wddPrty} and
equation~~\eqref{proport_eq}
$$
\leb\left(A_q(\sigma,v_0)\right),\,\leb\left(A'_q(\sigma,v_0)\right)\ge
\frac{1}{|R|}\frac{qd}{L}\frac{c}{q}\ge \frac{dc^2}{(2+c)|R|}.
$$
Thus, both $v_0$ and $v_0+\sigma$ are quasi-periods for the cocycle
$(\Phi_n)$. By Lemma \ref{lem-period}, they are periods. Hence
$\sigma$ is a period for the cocycle $(\Phi_n)$. Since
$\si\in\Sigma(\Phi)$ was arbitrary, the group of periods for
$\rho_{\alpha,\Phi}$ contains $\Si(\Phi)$.

\vskip 3mm If $H\subset G$ is a subgroup, we will denote by ``bar'' the reduction
modulo $H$. Thus, $\Phi:X\to G/H$ is a piecewise constant function.
The dynamical system $(X\times G/H,  \rho_{\alpha,\bPhi},\bmu)$ is the skew product
over $\rho_{\alpha}$ with the fiber $G/H$ and the displacement
function $\bPhi$. To simplify notation, we will denote it by $\overline{\rho_{\alpha,\Phi}}$.

Let $H\subset G$ (resp. $H'\subset G$) be the group
generated by $\Si(\Phi)$ (resp. $R(\Phi)$). We have shown that $H$
is contained in the group of periods for $\rho_{\alpha,\Phi}$. 
The function $\bPhi:X\to G/H$ is constant. Let
$a\in H'$ be such that $\bPhi=\bar{a}$. Then
$$
\overline{\rho_{\alpha,\Phi}}(x,\bg)=(\rho_{\al}(x),\bg+\bar{a}).
$$

Observe that $H'/H\subset G/H$ is the cyclic group generated by
$\bar{a}$. If $|H'/H|=\infty$, then $\rho_{\alpha,\bPhi}$ is
dissipative, contrary to Corollary~~\ref{rec1}. Thus, $H'/H$ is a
finite cyclic group. Let $|H'/H|=n$.

Let $f$ be a $\rho_{\alpha,\Phi}$-invariant measurable function.
Then $f$ defines a $\overline{\rho_{\alpha,\Phi}}$-invariant function; we
denote it by $\bar{f}$. By the above equation
$$
\bar{f}(\rho^n_{\al}(x),\bg)=\bar{f}(x,\bg).
$$
Hence $\bar{f}(x,\bg)$ depends only on $\bg$. In the self-explanatory
notation, $\bar{f}(x,\bg)=\bar{f}(\bg)$. Therefore, $\bar{a}$ is a
period for $\bar{f}$, and hence $a$ is a period for $f$.
\end{proof}
\end{thm}

The following is immediate from Proposition~~\ref{wdd-ri2} and
Theorem~~\ref{wdd-erg-thm}.

\begin{cor} \label{D-wdd}
Let $\Phi:X\to G$ be a piecewise constant, centered function such
that $R(\Phi)$ generates $G$. Let ${\mathcal D}\subset X$ be the set
discontinuities of $\Phi$. Let $0<\al<1$ be irrational. Suppose that
${\mathcal D}\subset\{\frac12(\Z\al+\Z) \setminus
(\Z\al+\Z)\}\mod1$. If $\alpha$ has property (D), then the skew
product $\rho_{\alpha,\Phi}$ is ergodic.
\end{cor}

\begin{rem}     \label{generic_rem}
{\em Let $\Psi:X\to\ZZ$ be the piecewise constant  function that
arose in our analysis of the rectangular Lorenz gas. See
equation~~\eqref{cocy_Psi}. It satisfies the  assumptions of
Corollary~~\ref{D-wdd}. Almost every irrational $\alpha$ has
property (D). Thus, Corollary~~\ref{D-wdd} implies the claim of
Theorem~~\ref{erg_decom_thm} for generic small obstacles. The
results in the next section will allow us to prove
Theorem~~\ref{erg_decom_thm} for arbitrary small obstacles.

}
\end{rem}

\subsection{Removing the genericity assumptions}      \label{spec_cocy_sub}
\hfill \break Let $0<\al<1$. Set
\begin{equation}  \label{gamma_zeta_eq}
\gamma = 1_{[0,{1\over 2}[} - 1_{[{1\over 2}, 1[},\ \zeta = 1_{[0,
{1\over 2} - {\alpha \over 2}[} - 1_{[{1\over 2}, 1- {\alpha \over
2}[}.
\end{equation}
Thus, $\ga,\zeta:X\to\Z$ are  piecewise constant functions. See Figure~~\ref{bill-fig2}.

\vspace {-4.cm}
\begin{figure} [htbp] 
\begin{center}
\includegraphics[scale=0.65]{bill-fig2.tex}
\end{center}
\vspace {1.cm}
\caption{The function $\zeta= 1_{[0, \beta[} - 1_{[{1 \over 2},
\beta + {1 \over 2}[}$ for $\be<1/2$.}   \label{bill-fig2}
\end{figure}

\vspace {1.cm}


We will establish the ergodicity of cocycles $(\ga_n),(\zeta_n)$ for
any irrational $\al$. First, we explain the heuristics. Suppose that
$\alpha$ does not satisfy property (D). Then, passing to a
subsequence, if need be, we have $a_n\to\infty$. Suppose that for
all sufficiently large $n$ the numbers $q_{2n}$ and $p_{2n}$ are
odd, while $q_{2n+1}$ is even and $p_{2n+1}$ is odd. Then for
$n>n_0$ the inequality $\zeta(q_n,\cdot)\ne 0$ holds only on sets of
very small measure. Thus, we cannot use the method of
Theorem~~\ref{wdd-erg-thm}. Instead, we will consider
$\zeta(tq,\cdot)$ for such $t$ that $||tq \alpha ||$ is close to zero
but big enough to ensure that $\zeta(tq,x) = \sum_{k=0}^{t-1}
\zeta(q,x+ kq\alpha)$ does not vanish on a set of measure bounded
away from zero. 

\newpage

We will informally refer to this idea as the ``filling method''.
Figure~~\ref{bill-fig6} illustrates it.

\vspace {-4 cm}
\begin{figure} [htbp]
\begin{center}
\includegraphics[scale=0.60]{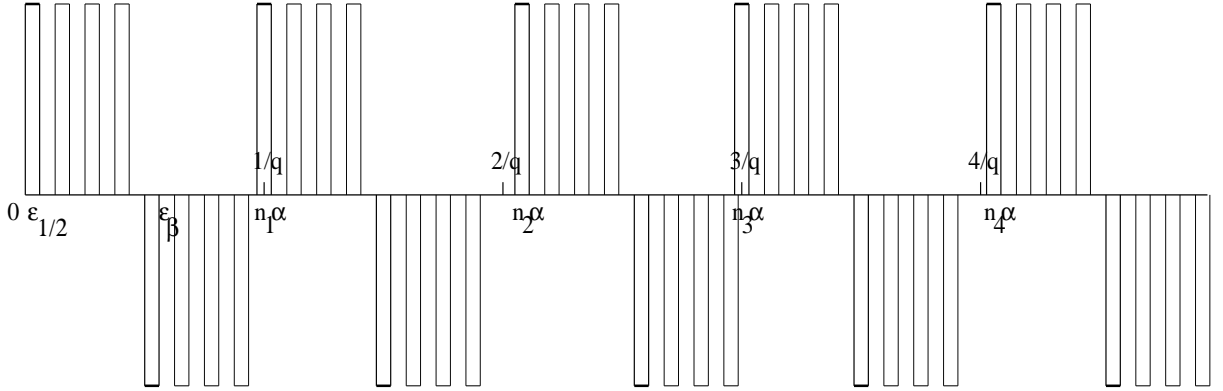}
\end{center}
\caption{\it Illustration of the filling method: Graph of
$\zeta(t_n q_n,\cdot)$.}                       
\label{bill-fig6}
\end{figure}

From now on, $0<\al<1$ is an arbitrary irrational number. We leave
the proof of the following lemma to the reader.

\begin{lem}    \label{undemi-lem}
If $q$ is odd and $q\|q\alpha\| < 1/2$, then for all $x$ we have
$\sum_{j=0}^{q-1}\ga(x +j\al)=\pm1$.
\end{lem}

By Lemma~~\ref{undemi-lem}, $1$ is a quasi-period for the cocycle
$(\gamma_n)$ over $\rho_{\al}$. By Lemma~~\ref{lem-period},
$(\gamma_n)$ is ergodic.


\begin{thm}      \label{erg-beta}
Let $\zeta:X\to\Z$ be given by equation~~\eqref{gamma_zeta_eq}. Then
the cocycle $(\zeta_n)$ over $\rho_{\al}$ is ergodic.
\begin{proof} Let $p_n/q_n$ be the convergents of $\alpha$. Let $p_n',q_n'\in\N$
be such that $q_n = 2q_n'$ or $q_n =2q_n'+1$ and $p_n = 2p_n'$ or
$p_n=2p_n'+1$, depending on the parity. Set $\alpha = {p_n/q_n} +
\theta_n$.

The set of discontinuities of $\zeta$ is $\{0,\beta=\frac12 -
\frac{\alpha}{2},\frac12,\beta'= 1 - \frac{\alpha}{2}\}$; the
respective jumps are $1,\,-1,\,-1,\,1$. If $t \in \{0, \beta,
\frac12, \beta' \}$, the corresponding discontinuities of
$\zeta_{q}$ are $\{ t - j \alpha: j= 0,..., q-1\}$.

Depending on the parities of $p_n,q_n$, we define partitions
$\{0, \beta, \frac12, \beta'\}=P_1\cup P_2$ as follows:

\medskip

\noindent 1) For $q_n$ odd and  $p_n$ even, we set $P_1 = \{0,
\beta' \}$, $P_2 = \{{1 \over 2}, \beta \}$;

\noindent 2) For $q_n$ even and $p_n$ odd, we set $P_1 = \{0,
{1\over 2} \}$, $P_2 = \{\beta, \beta'\}$;

\noindent 3) For $q_n$ odd and $p_n$ odd, we set $P_1 = \{0, \beta
\}$, $P_2 = \{{1\over 2}, \beta' \}$.

\medskip

Discontinuities of $\zeta_{q_n}$ which come from points in the same
atom of the partition are very close to each other; discontinuities
which come from points in distinct atoms of the partition are well
separated from each other.

We will consider in detail only case 2). The analysis of other cases
is similar, and we leave it to the reader. In what follows, all
numbers and equalities are understood ${\rm { \ mod \ 1}}$. To
simplify notation, we will suppress the subscript $n$. Thus, we
write $\alpha = p/q +\theta$, etc.

\medskip

\noindent a) The set of discontinuities, $D_0\subset X$,
corresponding to $t=0$ is $D_0=\{- j \alpha: j= 0, ..., q-1\}$. For
each integer $r$ there is $j_1(r) \in \{0, ..., q-1\}$ such that
$-j_1(r)p=r {\rm { \ mod \ }} q$. Hence, $D_0=\{ {r \over q} -
j_1(r) \theta: \ r= 0, ..., q-1 \}$.

\noindent b) Let $t=\frac12$. The corresponding set of
discontinuities is $D_{\frac12}=\{\frac{(q' -jp)}{q} - j\theta: j=
0, ..., q-1\}$. For each integer $r$ there is $j_2(r) \in \{0, ...,
q-1\}$ such that $q'-j_2(r)p=r {\rm { \ mod \ }} q$. Thus,
$D_{\frac12}=\{ {r \over q} - j_2(r) \theta, \ r= 0, ..., q-1 \}$.

\noindent c) Let $t=1 - \frac{\alpha}{2}$. Then the set of
discontinuities is $D_{1 - \frac{\alpha}{2}}=\{\frac{1}{2q} -
\frac{(p'+1 +jp)}{q} - (j+ \frac12) \theta: j= 0, ..., q-1\}$. For
each integer $r$ there is $j_3(r) \in \{0, ..., q-1\}$ such that
$-(p'+1+j_3(r)p)=r {\rm \ mod \ } q$. Hence $D_{1 -
\frac{\alpha}{2}}=\{ \frac{r}{q} +  \frac{1}{2q} - (j_3(r) +
\frac12) \theta: \ r= 0, ..., q-1 \}$.

\noindent d) Let $t=\frac12 - \frac{\alpha}{2}$. The set of
discontinuities is $D_{\frac12 - \frac{\alpha}{2}}=\{\frac{1}{2q} -
\frac{(-q'+p'+1 + jp)}{q} - (j + \frac12) \theta: j= 0,\dots,q-1\}$.
For each integer $r$ there is $j_4(r) \in \{0, ..., q-1\}$ such that
$q'-(p'+1+j_4(r)p)=r {\rm \ mod \ } q$. Hence $D_{\frac12 -
\frac{\alpha}{2}}=\{\frac{r}{q} + \frac{1}{2q} - (j_4(r) + \frac12)
\theta: \ r= 0, ..., q-1 \}$.

The set of discontinuities of $\zeta_{q}$ is $D_0\cup
D_{\frac12}\cup D_{1 - \frac{\alpha}{2}}\cup D_{\frac12 -
\frac{\alpha}{2}}$. Observe that in all cases $|j_i(r) \theta| \leq
|q\theta|$; since $(j_2 - j_1) p = \frac12 q {\rm \ mod \ } q$ and
$(j_4 - j_3) p = \frac12 q {\rm \ mod \ } q$, we have
\begin{eqnarray}
j_2(r) = j_1(r) \pm {1\over 2} q, \ j_4(r) = j_3(r) \pm
{1\over 2} q. \label{j1j2}
\end{eqnarray}

We are going to determine the values taken by the cocycle
$\zeta_q(x)$ for $x$ in a neighborhood of the typical interval
$[\frac{r}{q}, \frac{r+1}{q}]$, where $r$ is an integer in $\{1,
..., q-1\}$. Assume, for concreteness, that $\theta < 0$, $j_1 =
j_2+ \frac12 q$, $j_4 = j_3 + \frac12 q$.\footnote{The analysis for
$\theta >0$ and/or $j_1 = j_2 - \frac12 q$, and/or $j_4 = j_3 -
\frac12 q$ is analogous.} Let $x$ start at $\frac{r}{q}$ and let it
move to the right; set $\zeta_q(x)=a$. The value of the cocycle
$\zeta_q(x)$ is constant until $x$ crosses the discontinuity
(corresponding to $t=0$) at $\frac{r}{q} - j_1(r)\theta$, where
$\zeta_q(x)$ increases by $1$. After that the cocycle does not
change until $x$ crosses the discontinuity at $\frac{r}{q} - j_2(r)
\theta$ (corresponding to $t=\frac12$) where the cocycle decreases
by $1$, returning to the value $a$.

The first two discontinuities occur before $x$ crosses
the two other discontinuities under the condition that $|j_i(r)
\theta|$ is less than ${1 \over 2q}$. This takes place if
$q^2|\theta| < \frac12$, a condition which holds below because we
consider the case when $q^2|\theta|$ is small. As $x$ continues to
move to the right, the cocycle remains at the value $a$ until, near
$\frac{r}{q} + \frac{1}{2q}$, it increases by $1$ at the point
$\frac{r}{q} + \frac{1}{2q} - (j_3(r) + \frac12) \theta$, a
discontinuity corresponding to $t=1 - \frac{\alpha}{2}$, and then
decreases by $1$ at $\frac{r}{q} + \frac{1}{2q} - (j_4(r) + \frac12)
\theta$, a discontinuity corresponding to $t=\frac12 -
\frac{\alpha}{2}$.

Therefore, we have
\begin{eqnarray*} \zeta_{q}&=&a \pm 1 {\rm \ on \ }
]\frac{r}{q} - j_1(r) \theta,\frac{r}{q} - j_2(r) \theta[, \\
\zeta_{q}&=& a \pm 1 {\rm \ on \ } ]\frac{r}{q} + \frac{1}{2q} -
(j_3(r) + \frac12) \theta, \frac{r}{q} + \frac{1}{2q} - (j_4(r) +
\frac12) \theta[.
\end{eqnarray*}
Elsewhere, $\zeta_{q}= a$.

This analysis is valid for every interval $[\frac{r}{q},
\frac{r+1}{q}]$. The order of the discontinuity points may change,
but not the order between the groups of discontinuity. In
particular, this implies that $\zeta_{q}=a$ on a subset of large
measure in $[0,1]$. Since the mean value of $\zeta_{q}$ is zero, we
have $a=0$.

\medskip 

We will now finish the proof. If $\alpha$ satisfies condition (D),
then the claim holds, by Corollary~~\ref{D-wdd}. Thus, we assume
that $\alpha$ does not satisfy condition (D). Then one of the cases
1), 2), 3) materializes for an infinite sequence $(n_k)$ such that
$a_{n_k+1} \to \infty$. We will then say, for brevity, that a case
{\em occurs infinitely often}. If case 1) occurs infinitely often,
then $1$ is a quasi-period for $(\zeta_n)$.  See figure~~\ref{bill-fig3}.

\vspace{0.9 cm}
\begin{figure} [htbp]
\begin{center}
\includegraphics[scale=0.50]{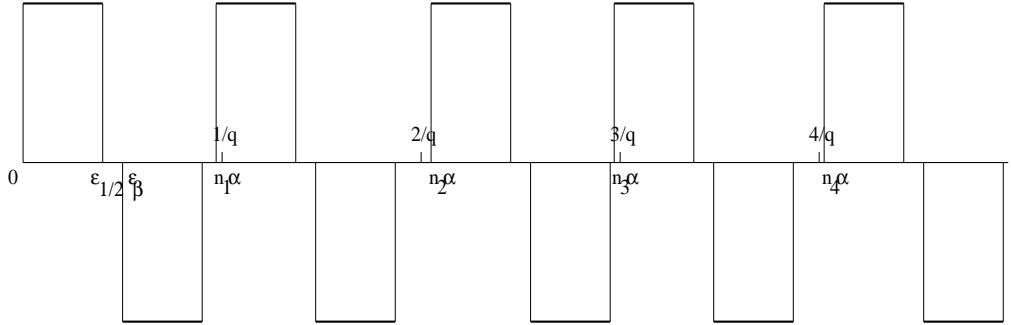}
\end{center}
\caption{\it Graph of the cocycle $\zeta_q$ for $q$
odd, $p$ even. Its value is $1$ on a set of measure $\geq \delta
> 0$.}  
\label{bill-fig3}
\end{figure}

Suppose now that case 2) occurs infinitely often. We will use the
preceding analysis. For $1\le r \le q_{n_k}-1$ set
\begin{eqnarray*}
I_{k,r} &=& ]\frac{r}{q_{n_k}} - j_1(r) \theta_{n_k},
\frac{r}{q_{n_k}} - j_2(r) \theta_{n_k}[,\\
J_{k,r} &=&  ]\frac{r}{q_{n_k}} + \frac{1}{2q_{n_k}} - (j_3(r) +
\frac12) \theta_{n_k}, \frac{r}{q_{n_k}} + \frac{1}{2q_{n_k}} -
(j_4(r) + \frac12) \theta_{n_k}[.
\end{eqnarray*}

By equation~~\eqref{j1j2}, these intervals have length $\frac12 q_{n_k}
|\theta_{n_k}|$. At the scale $\frac{1}{q_{n_k}}$, they are close to
$\frac{r}{q_{n_k}}$ and to $\frac{r}{q_{n_k}} + \frac{1}{2q_{n_k}}$
respectively. Outside of these intervals, $\zeta(q_{n_k}, \cdot)=0$.

\vskip 3mm Let $\delta \in ]0, {1\over 4}[$. Set $t_{k} = \lfloor
\delta a_{{n_k}+1}\rfloor$. For $J\subset X$ and $u\in\R$, we set
$(J + u)=J+u \mod 1$. Set
$$A_k = \cup_{j=0}^{q_{n_k}-1} \, \cup_{s=0}^{t_k-1} \, (I_{k,j} - sq_{n_k}
\alpha), \ B_k =  \cup_{j=0}^{q_{n_k}-1} \, \cup_{s=0}^{t_k-1} \,
(J_{{k},j} - sq_{n_k}\alpha).$$

The distance between the intervals $I_{{k},r}$ and $J_{{k},r}$ is at
least $\frac{1}{2q_{n_k}} - q_{n_k}|\theta_{n_k}|$. Since, by the
choice of $t_k$, we have $q_{n_k} |\theta_{n_k}| t_k \leq
\frac{1}{2q_{n_k}} - q_{n_k}|\theta_{n_k}|$, the translated
intervals in the definition of $A_k$ and $B_k$ do not overlap.

\vskip 3mm Let us consider the cocycle at time $t_{k}q_{n_k}$. We
have $\zeta(t_{k}q_{n_k}, x) = \sum_{s=0}^{t_k-1} \zeta(q_{n_k}, x +
s q_{n_k}\alpha)$. By the preceding analysis of the values of
$\zeta_q$, we have $\zeta(t_k q_{n_k}, \cdot)= \pm 1$ on $A_k$ and
$B_k$. Also
$$
\leb(A_k) = {1\over 2} t_k q_{n_k} q_{n_k}|\theta_{n_k}| \geq
{1\over 2} \delta a_{{n_k}+1} \frac{q_{n_k}}{q_{{n_k}+1}} \geq
{1\over 2} \delta > 0.
$$

Since $t_k q_{n_k} \alpha {\rm \ mod \ } 1 \to 0$, and since on
$A_k$ the cocycle takes at most the two values $\pm 1$, we have
shown that $1$ or $-1$ is a quasi-period for the cocycle
$(\zeta_n)$.

\vskip 3mm The possibility that case 3) occurs infinitely often is
analyzed the same way, with the conclusion that $1$ is a
quasi-period.

Thus, no matter which of the three cases occurs infinitely often, $1$
is a quasi-period for the cocycle $(\zeta_n)$. The claim now
follows, by Lemma~~\ref{lem-period}.
\end{proof}
\end{thm}

\vskip 3mm We will now establish the main result of this subsection.

\begin{thm} \label{ergo-Psi}
Let $\Psi:X\to\ZZ$ be the function defined by
equation~~\eqref{cocy_Psi}. Then, over any irrational rotation
$\rho_{\al}: x \to x+\alpha {\rm \ mod \ } 1$ on the circle, the
corresponding cocycle is ergodic.
\begin{proof}
Set $\Psi=(\psi_1,\psi_2)$ and $\be= \frac12 - \frac{\al}{2}$. The
functions $\psi_1, \psi_2$ satisfy
$$\psi_1(x)+\psi_2(x) = \gamma(x),\ \psi_1(x) -\psi_2(x) =
\gamma(x+\beta).$$

By Lemma~~\ref{4odd} and Lemma~~\ref{undemi-lem}, for an infinite
sequence $(n_k)$ we have $\gamma(q_{n_k},x), \gamma(q_{n_k},
x+\beta)\in\{\pm 1\}$. This implies the existence of measurable sets
$A_{k} \subset X$ satisfying $\leb(A_{k})> \frac14$, and such that
for $x \in A_{k}$ the vector function $(\gamma(q_{n_k},x),
\gamma(q_{n_k},x+\beta))$ is constant. Its values are $(+1,+1)$,
($+1, -1)$, $(-1,+1)$, or $(-1,-1)$.

Thus, for $x\in A_{k}$, the vector function $(\psi_1(q_{n_k},x), \psi_2(q_{n_k},x))$ is
identically $(1,0)$, $(0,1)$, $(-1,0)$, or $(0,-1)$. Hence, one of
the elements $(\pm 1,0),(0,\pm1)\in\ZZ$ is a quasi-period for the
cocycle $(\Psi_n)$. Suppose, for instance, that $(1,0)$ is a
quasi-period. Hence, by Lemma~~\ref{lem-period}, $(1,0)$ is a
period. Let $f$ be a $\rho_{\alpha, \Psi}$-invariant function. It
defines a $\rho_{\alpha, \zeta}$-invariant function on $X \times\Z$.
Here $\zeta:X\to\Z$ is given by equation~~\eqref{gamma_zeta_eq}. By
Theorem~~\ref{erg-beta}, $f=\const$, i. e.,  the cocycle $(\Psi_n)$
is ergodic. The other cases are disposed of the same way.
\end{proof}
\end{thm}

\medskip

\noindent{\bf Acknowledgements} In the course of preparation of this
work the coauthors have made several visits to each other's home
institutions. It is a pleasure to thank IRMAR and UMK for making
these visits possible. The work of E.G. was partially supported
by MNiSzW grant N N201 384834.

\end{document}